\documentclass{amsart}
\usepackage{amsmath,amscd,amssymb,latexsym,graphicx}
\usepackage{enumerate}

\usepackage{color}

\newtheorem{theorem}{Theorem}[section]

\newtheorem{dfprop}[theorem]{Definition-Proposition}

\newtheorem{proposition}[theorem]{Proposition}
\newtheorem{corollary}[theorem]{Corollary}

\theoremstyle{definition}
\newtheorem{definition}[theorem]{Definition}
\newtheorem{remark}[theorem]{Remark}
\newtheorem{remarks}[theorem]{Remarks}

\newtheorem{example}[theorem]{Example}
\newtheorem{examples}[theorem]{Examples}

\numberwithin{equation}{section}

\reversemarginpar

\newcommand{\cA}{\mathcal{A}}
\newcommand{\cB}{\mathcal{B}
}\newcommand{\cC}{\mathcal{C}}

\newcommand{\cE}{\mathcal{E}}
\newcommand{\cF}{\mathcal{F}}
\newcommand{\cG}{\mathcal{G}}
\newcommand{\cH}{\mathcal{H}}
\newcommand{\cI}{\mathcal{I}}
\newcommand{\cJ}{\mathcal{J}}
\newcommand{\cK}{\mathcal{K}}
\newcommand{\cL}{\mathcal{L}}
\newcommand{\cM}{\mathcal{M}}
\newcommand{\cN}{\mathcal{N}}

\newcommand{\cS}{\mathcal{S}}
\newcommand{\cT}{\mathcal{T}}

\newcommand{\cV}{\mathcal{V}}
\newcommand{\cW}{\mathcal{W}}

\newcommand{\fS}{ \mathsf{S}}

\newcommand{\NN}{\mathbb N}
\newcommand{\ZZ}{\mathbb Z}

\newcommand{\RR}{\mathbb R}
\newcommand{\CC}{\mathbb C}

\newcommand{\ind}{\mathop{\mathrm{Ind}}\nolimits}

\newcommand\Gr[1]{\cG^{(#1)}}
\newcommand{\pb}[1]{{}^*#1^*}
\newcommand{\To}{\longrightarrow}

\newcommand\smx{X^{\circ}}
\newcommand\smy{Y^{\circ}}
\newcommand\smE{E^{\circ}}
\newcommand\smf{f^{\circ}}
\newcommand\pa{\partial}
\newcommand{\ie}{{\em i.e.}\ }

\title{Groupoids and an index theorem for conical
pseudo-manifolds}

\author{Claire Debord}
\address{Laboratoire de Math{\'e}matique , Universit{\'e} Blaise Pascal,
  Campus universitaire des c{\'e}zeaux, 63177 Aubi{\`e}re Cedex, France}
\email{debord@math.univ-bpclermont.fr}
\author{Jean-Marie Lescure}
\address{Laboratoire de Math{\'e}matique , Universit{\'e} Blaise Pascal,
  Campus universitaire des c{\'e}zeaux, 63177 Aubi{\`e}re Cedex, France}
\email{lescure@math.univ-bpclermont.fr}
\author{Victor Nistor}
\address{Institutul de Matematic\u a al Academiei
Rom\^ane and Department of Mathematics, Pennsylvania State
University, University Park, PA 16802, USA}
\email{nistor@math.psu.edu}

\date{\today}

\begin{document}

\maketitle

\begin{abstract}
We define an analytical index map and a topological index map for
conical pseudomanifolds. These constructions generalize the analogous
constructions used by Atiyah and Singer in the proof of their
topological index theorem for a smooth, compact manifold $M$. A main
new ingredient in our proof is a non-commutative algebra that plays in
our setting the 
role of $\mathcal{C}_0(T^*M)$. We prove a Thom isomorphism
between non-commutative algebras which gives a new example of wrong
way functoriality in $K$-theory. We then give a new 
proof of the Atiyah-Singer index theorem using deformation groupoids
and show how it generalizes to conical pseudomanifolds. We thus prove
a topological index theorem for conical pseudomanifolds.
\end{abstract}

\tableofcontents

%%%%%%%%%%%%%%%%%%%%%%%%%%%%%%%%%%%%%%%%%%%%%%%%%%%%%%%%%%%%%%%%%%%%%%%%%%%%%%%%%%%%%%%%

\section*{Introduction}

Let $V$ be a closed, smooth manifold and let $P$ be an elliptic
pseudo-differential operator acting between Sobolev spaces of
sections of two
vector bundles over $V$. The ellipticity of $P$ ensures that $P$
has finite dimensional kernel and cokernel. The difference
\begin{equation*}
    \ind P:=\mbox{dim} (Ker P) -\mbox{dim} (Coker P)
\end{equation*}
is called the {\em Fredholm index} of $P$ and turns out to
depend only on the $K$-theory class $[\sigma(P)]\in K^0(T^* V)$ of
the principal symbol of $P$ (we always use $K$-theory with
compact supports). Since every element in $ K^0(T^* V)$ can be
represented by the principal symbol of an elliptic
pseudo-differential operator, one obtains in this way a group morphism
%\begin{equation*}
%\begin{array}{ccc}
%    K^0(T^*V) & \longrightarrow & \ZZ \\
%    \lbrack \sigma(P) \rbrack & \mapsto & \ind P
%\end{array}
%\end{equation*} 
\begin{equation}
        \ind_a^V :  K^0(T^*V)  \longrightarrow  \ZZ , \quad 
        \ind_a^V(\sigma(P)) =  \ind P,
        \end{equation}
called the {\it {analytical index}}, first introduced by M. Atiyah
and I. Singer \cite{AS1}. 

At first sight, the map $\ind_a^V$ seems to depend essentially on the analysis 
of elliptic equations, but the main result in \cite{AS1} is that the index map can
be also defined in a purely topological terms in terms of embedding of
$V$ in Euclidean spaces. This definition in terms on embeddings leads to the so 
called {\it topological index} of Atiyah and Singer and the main result of
\cite{AS1} is that the topological index map $\ind_t^V$ and the Fredholm index
map $\ind_a^V$ coincide. See \cite{Carvalho} for review of these results,
including an extension to non-compact manifolds.

The equality of the topological and Fredholm indices then allowed
M. Atiyah and I. Singer to obtain a formula for the index of an
elliptic operator $P$ in terms of the Chern classes of
$[\sigma(P)]$. Their formula, the celebrated Atiyah-Singer Index
Formula, involves, in addition to the Chern character of the principal
symbol of $P$, also a universal characteristic class associated with
the manifold, the so called Todd class of the given manifold.

It is a natural and important question then to search for extensions
of the Atiyah--Singer results. It is not the place here to mention all
existing generalizations of the Atiyah-Singer index theory, but let us
mention here the fundamental work on Connes on foliations
\cite{Connes1, Connes2, CoBook, ConnesMoscovici, CS} as well as
\cite{BenameurNistor, Gorokhovsky, NistorFol, NistorDoc}. The index
theorem for families and Bismut's superconnection formalism play an
important role in the study of the so called ``anomalies'' in physics
\cite{Bismut1, Bismut2, BismutFreed, FreedWitten}. A different but
related direction is to extend this theory to singular spaces
\cite{APS1}. An important step in the index problem on singular
manifolds was made by Melrose \cite{meaps, MelroseScattering} and
Schulze \cite{Schu1,Schu2} who have introduced the ``right class of
pseudodifferential operators'' for index theory on singular
spaces. See also \cite{aln2, FedSchu, FedSchuTark, Schrohe, Schulze,
SchuSterSha}.  Generalizations of this theory to singular spaces may
turn out to be useful in the development of efficient numerical
methods \cite{BNZ1}.

In this paper, we shall focus on the case of a pseudomanifold $X$ with
isolated conical singularities. In earlier work \cite{DL}, the first
two authors defined a $C^*$-algebra $A_X$ that is dual to the algebra
of continuous functions on $X$ from the point of view of $K$-theory
(\ie $A_X$ is a ``$K$-dual of $X$'' in the sense of
\cite{CoBook,CS,Ka2}), which implies that there exists a natural
isomorphism
\begin{equation}
\label{symbolmap}
    K_0(X)\overset{\Sigma_X}\longrightarrow K_0( A_X )
\end{equation}
between the $K$-homology of $X$ and the $K$-theory of $T^{\fS}X$.  The
$C^*$-algebra $A_X$ is the $C^*$-algebra of a groupoid denoted
$T^{\fS}X$.

One of the main results in \cite{jml}, see also
\cite{MelPiaz,NSS1,NSS2,Sav1} for similar results using different
methods, is that the inverse of the map $\Sigma_X$ of Equation
\eqref{symbolmap} can be realized, as in the smooth case, by a map
that assigns to each element in $K_0(A_X)$ an elliptic operator. Thus
elements of $K_0(A_X)$ can be viewed as the symbols of some natural
elliptic pseudodifferential operators realizing the $K$-homology of
$X$. Of course, in the singular setting, one has to explain what is
meant by ``elliptic operator'' and by ``symbol'' on $X$. An example of
a convenient choice of elliptic operator in our situation is an
elliptic pseudodifferential operator in the $b$-calculus
\cite{meaps,Schu1} or Melrose's $c$-calculus. As for the symbols, the
notion is more or less the same as in the smooth case. On a manifold
$V$, a symbol is a function on $T^*V$. For us, it will be convenient
to view a symbol as a pointwise multiplication operator on
$C_c^{\infty}(T^* V)$. A Fourier transform will allow us then to see a
symbol as a family of convolution operators on $C_c^{\infty}(T_x V),\
x\in V$. Thus symbols on $V$ appear to be pseudo-differential
operators on the {\it groupoid} $TV$. This picture generalizes then
right away to our singular setting. In particular, it leads to a good
notion of symbol for conical pseudomanifolds and enables us to
interpret \ref{symbolmap} as the principal symbol map.

In order to better explain our results, we need to introduce some
notation. If $G$ is an amenable groupoid, we let $K^0(G)$ denote
$K_0(C^*(G))$. The analytical index is then defined exactly as in the
regular case by
\begin{equation*}\begin{array}{cccc}
    \ind_a^{X} : & K^0(T^{\fS}X) &
    \rightarrow & \ZZ \\ & [a] &
    \mapsto & \ind(\Sigma^{-1}_X(a)),
\end{array}\end{equation*} 
where $\ind:K_0(X)\rightarrow \ZZ$ is the usual Fredholm index
on compact spaces. Moreover one can generalise the {\it tangent
groupoid} of A. Connes to our situation and get a nice description of
the analytical index.

Following the spirit of \cite{AS1}, we define in this article a
topological index $\ind_t^{X}$ that generalizes the classical
one and which satisfies the equality:
\begin{equation*}
 \ind_a^{X}=\ind_t^{X}.
\end{equation*} 
In fact, we shall see that all ingredients of the classical
topological index have a natural generalisation to the singular
setting.\smallskip

\noindent $\bullet$ Firstly the embedding of a smooth manifold into
$\RR^N$ gives rise to a normal bundle $N$ and a Thom isomorphism
$K^0(T^*V)\rightarrow K^0(T^*N)$. In the singular setting we embed $X$
into $\RR^{N}$, viewed as the cone over $\RR^{N-1}$. This gives rise
to a {\it conical vector bundle} which is a conical pseudomanifold
called the {\it normal space} and we get an isomorphism:
$K^0(T^{\fS}X)\rightarrow K^0(T^{\fS}N)$. This map restrict to the
usual Thom isomorphism on the regular part and is called again the
{\it Thom isomorphism}.

\noindent $\bullet$ Secondly, in the smooth case, the normal bundle
$N$ identifies with an open subset of $\RR^N$, and thus provides an
excision map $K(TN)\rightarrow K(T\RR^N)$. The same is true in the
singular setting : $T^{\fS}N$ appears to be an open subgroupoid of
$T^{\fS}\RR^N$ so we have an excision map $K^0(T^{\fS}N)\rightarrow
K^0(T^{\fS}\RR^N)$.

\noindent $\bullet$ Finally, using the Bott periodicity
$K^0(T^*\RR^N)\simeq K^0(\RR^{2N}) \rightarrow \ZZ$ and a natural
$KK$-equivalence between $T^{\fS}\RR^N$ with $T\RR^N$ we obtain an
isomorphism $K^0(T^{\fS}\RR^N)\rightarrow \ZZ$.\smallskip

\noindent As for the usual definition of the topological index,
this allows us to define our generalisation of the topological $\ind_t$
for conical manifolds.

This construction of the topological index is inspired from
the techniques of {\it deformation groupoids} introduced by M. Hilsum
and G. Skandalis in \cite{HS1}. Moreover, the demonstration of the
equality between $\ind_a$ and $\ind_t$ will be the same in
the smooth and in the singular setting with the help of deformation
groupoids.

We claim that our index maps are straight generalisations of the
classical ones. To make this claim more concrete, consider a closed
smooth manifold $V$ and choose a point $c\in V$. Take a neighborhood
of $c$ diffeomorphic to the unit ball in $\RR^n$ and consider it as
the cone over $S^{n-1}$. This provides $V$ with the structure of a
conical manifold. Then the index maps $\ind^{\fS}_*:K^0(
T^{\fS}V)\rightarrow \ZZ$ and $\ind_*:K^0( TV)\rightarrow \ZZ$
both correspond to the canonical map $K_0(V)\rightarrow \ZZ$ through
the Poincar{\'e} duality $K_0(V)\simeq K^0(T^*V)$ and $K_0(V)\simeq
K^0(T^{\fS}V)$. In other words both notions of indices coincide trough
the $KK$-equivalence $TV\simeq T^{\fS}V$.

We will investigate the case of general stratifications and
the proof of an index formula in forthcoming papers.
 
The paper is organized as follows.  In Section \ref{geom} we describe
the notion of conical pseudomanifolds and conical bundles. Section
\ref{sec2} reviews general facts about Lie groupoids. Section
\ref{sec3} is devoted to the construction of tangent spaces and
tangent groupoids associated to conical pseudomanifolds as well as
other deformation groupoids needed in the subsequent sections.
Sections \ref{sec4} and \ref{Thominv} contain the construction of
analytical and topological indices, and the last section is devoted to
the proof of our topological index theorem for conical
pseudomanifolds, that is, the proof of the equality of analytical and
topological indices for conical pseudomanifolds.

%%%%%%%%%%%%%%%%%%%%%%%%%%%%%%%%%%%%%%%%%%%%%%%%%%%%%%%%%%%%%%%%%%%%%%%%%%%%%%%%%%%%%%%%

%\section{Basic definitions}

\section{Cones and stratified bundles\label{geom}}

We are interested in studying conical pseudomanifolds, which are
special examples of stratified pseudomanifolds of depth one
\cite{GoMa}. We will use the notations and equivalent descriptions
given by A. Verona in \cite{Ve} or used by J.P. Brasselet, G.
Hector and M. Saralegi in \cite{BHS}. See \cite{HW} for a review
of the subject.

\subsection{Conical pseudomanifolds} If $L$ is a smooth manifold,
the {\it cone} over $L$ is, by definition, the topological space
\begin{equation}\label{eq.def.cone}
    cL:=L\times [0,+\infty [/ L\times\{0\}.
\end{equation}
Thus $L\times\{0\}$ maps into a single point $c$ of $cL$. We shall
refer to $c$ as the {\it singular point of $L$}. If $z\in L$ and
$t\in [0,+\infty[$ then $[z,t]$ will denote the image of $(z,t)$
in $cL$. We shall denote by
\begin{equation*}
    \rho_{cL}:cL \rightarrow [0,+\infty[, \quad \rho_{cL}([z, t])
    := t
\end{equation*}
the map induced by the second projection and we call it the {\it
defining function} of the cone.

\begin{definition}
A {\it conical stratification} is a triplet $(X,\fS,\cC)$ where
\begin{enumerate}[\rm (i)]
\item $X$ is a Hausdorff, locally compact, and secound countable
space.

\item $\fS \subset X$ is a finite set of points, called
the {\it singular set} of $X$, such that $\smx :=X\setminus \fS$
is a smooth manifold.

\item $\cC = \{ (\cN_{s},\rho_s,L_s) \}_{s\in \fS}$ is the set
of {\it control data}, where $\cN_{s}$ is an open neighborhood of
$s$ in $X$ and $\rho_s: \cN_s \rightarrow [0,+\infty[$ is a
surjective continuous map such that $\rho_s^{-1}(0)={s}$.

\item For each $s \in \fS$, there exists a homeomorphism
$\varphi_s:\cN_s \rightarrow cL_s$, called {\it trivialisation
map}, such that $\rho_{cL_s} \circ \varphi_s  = \rho_s$ and such
that the induced map $\cN_s \setminus \{s\} \to L_s\times
]0,+\infty[$ is a diffeomorphism. Moreover, if $s_0,\ s_1 \in \fS$
then either $\cN_{s_0}\cap\cN_{s_1}=\emptyset$ or $s_0=s_1$.
\end{enumerate}
\end{definition}

\smallskip Let us notice that it follows from the definition that the
connected components of $\smx$ are smooth manifolds. These
connected components are called the {\it regular strata} of $X$.

\begin{definition} Two conical stratifications $(X,\fS_X,\cC_X)$
and $(Y,\fS_Y,\cC_Y)$ are called {\em isomorphic} if there is an
homeomorphism $f:X\rightarrow Y$ such that:
\begin{enumerate}[\rm (i)]
    \item $f$ maps $\fS_X$ onto $\fS_Y$,
    \item $f$ restricts to a smooth diffeomorphism $\smf : \smx
    \rightarrow \smy$,
    \item the defining function $\rho_s$ of any $s\in \fS_X$ is
    equal to $\rho_{f(s)}\circ f$, where $\rho_{f(s)}$ is the
    defining function of $f(s)\in \fS_Y$ (in particular
    $f(\cN_s)=\cN_{f(s)}$).
\end{enumerate}
An isomorphism class of conical stratifications will be called a
{\em conical pseudomanifold}.
\end{definition}

\smallskip In other words, a conical pseudomanifold is a locally compact,
metrizable, second countable space $X$ together with a finite set
of points $\fS \subset X$ such that $\smx = X \setminus \fS$ is a
smooth manifold and one can find a set of control data $\cC$ such
that $(X,\fS,\cC)$ is a conical stratification.

\smallskip Let $M$ be a smooth manifold with boundary $L := \pa M$. An easy
way to construct a conical pseudomanifold is to glue to $M$  the
{\it closed cone} $\overline{cL}:=L\times [0,1]/ L\times\{0\}$
along the boundary.

%\medskip \centerline{\includegraphics[width=7cm]{cone.pdf}}
\medskip \centerline{\includegraphics[width=7cm]{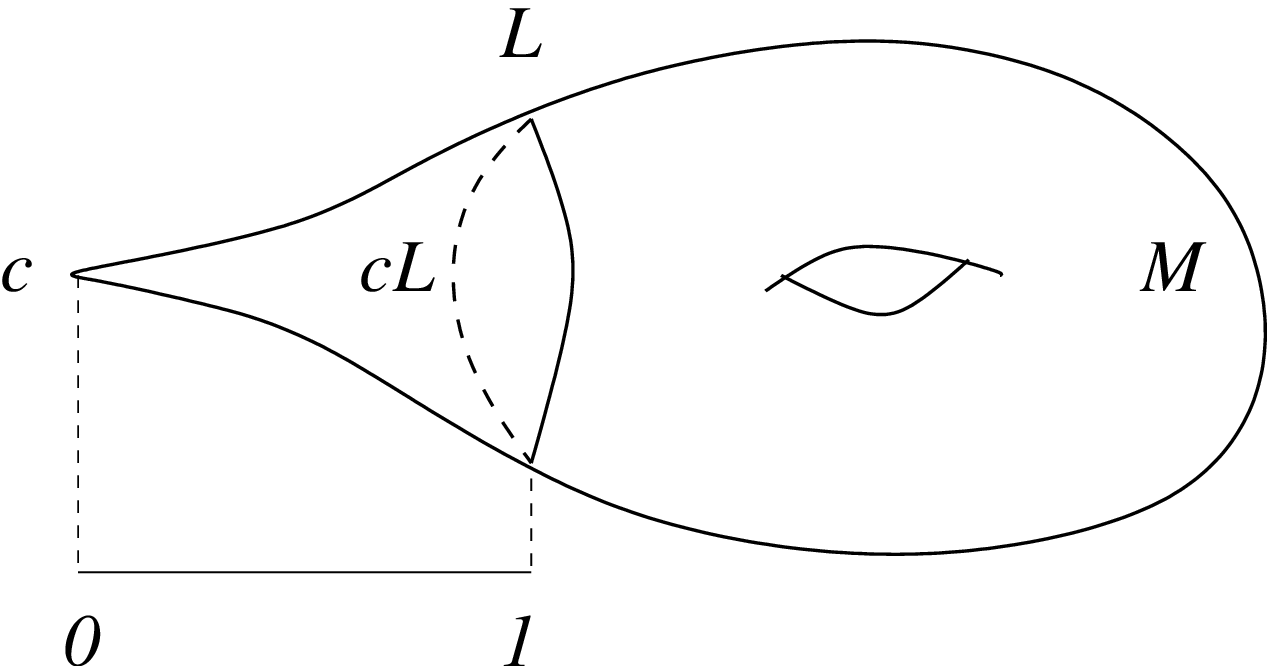}}

\smallskip \noindent Notice that we do not ask the link
$L$ to be connected. For example, if $M$ is a smooth manifold, the
space $M\times S^1 / M\times \{p\}$, $p \in S^1$, is a conical
pseudomanifold with $L$ consisting of two disjoint copies of $M$:

%\medskip \centerline{\includegraphics[width=5cm]{pinchedtorus.pdf}}
\medskip \centerline{\includegraphics[width=5cm]{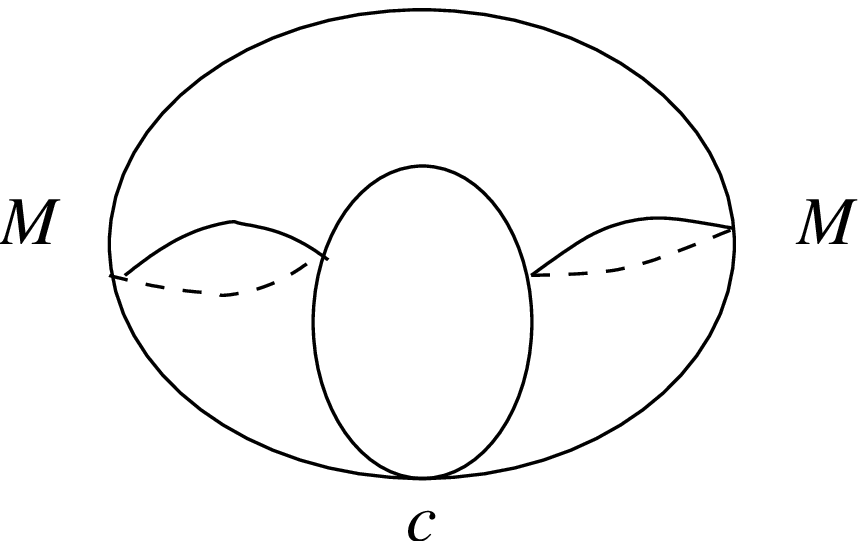}}

\subsection{Conical bundles}
We next introduce ``conical bundles,'' a class of spaces not to be
confused with vector bundles over conical manifolds. Assume that
$L$ is a smooth manifold, $cL$ is the cone over $L$,
$\pi_{\xi}:\xi \rightarrow L$ is a smooth vector bundle over $L$,
and $c\xi$ is the cone over $\xi$. We define $\pi: c\xi
\rightarrow cL$ by $\pi([z,t])=[\pi_{\xi}(z),t]$ for $(z,t)\in
\xi\times [0,+\infty[$. The set $(c\xi,\pi)$ is the {\it cone}
over the vector bundle $(\xi,\pi_{\xi})$. Let us notice that the
fiber above the singular point of $cL$ is the singular point of
$c\xi$. In particular, $c\xi$ is {\em not} a vector bundle over
$cL$.

\begin{definition}\label{def.con.bdle}
Let $(X,\fS_X,\cC_X)$ be a conical stratification. A {\it conical
vector bundle} $(E,\pi)$ over $X$ is a conical stratification
$(E,\fS_E,\cC_E)$ together with a continuous surjective map $\pi:E
\rightarrow X$ such that:
\begin{enumerate}
\item $\pi$ induces a bijection between the singular sets $\fS_E$
and $\fS_X$.
\item If $\smE :=E \setminus \fS_E$, the restriction
$\pi^{\circ}: \smE \rightarrow \smx$ is a smooth vector bundle.
\item The control data $\{\cM_{z},\rho_z,\xi_z\}_{z\in \fS_E}$ of
$E$ and $\{\cN_{s},\rho_s,L_s\}_{s\in \fS_X}$ of $X$ satisfies:
$\cM_{z}=\pi^{-1}(\cN_{\pi(z)})$ and $\rho_z=\rho_{\pi(z)}\circ
\pi$. Moreover for $z\in \fS_E$ and $s=\pi(z) \in \fS_X$, the
restriction $\pi_z: \cM_{z} \rightarrow \cN_{s} $ is a cone over
the vector bundle $\xi_z$. More precisely, we have the following
commutative diagram
\begin{equation*}
\begin{CD} \cM_{z} @>{\pi_z}>> \cN_{s} \\
@V{\Psi_z}VV @VV{\varphi_s}V \\ c\xi_z  @>>> cL_s \end{CD}
\end{equation*}
where $\xi_z\rightarrow L_s$ is a smooth vector bundle over $L_s$,
the bottom horizontal arrow is the cone over $\xi_z\rightarrow
L_s$  and $\Psi_z$, $\varphi_s$ are trivialisation maps.
\end{enumerate}
If $X$ is a conical pseudomanifold, the isomorphism class of a
conical vector bundle over a conical stratification
$(X,\fS_X,\cC_X)$ will be call again a conical vector bundle over
$X$.
\end{definition}

\smallskip We are interested in conical vector bundles because they
allow us to introduce the right notion of tubular neighborhood in
the class of conical manifolds.

\smallskip Let $L$ be a compact manifold and $cL$ the cone over
$L$. For $N \in \NN$ 
large enough, we can find an embedding $j_L : L\rightarrow S^{N-1}$ where
$S^{N-1}\subset \RR^{N}$ denotes the unit sphere.
Let $\cV_L \rightarrow L$ be the normal bundle of this embedding.
We let $c\cV_L=\cV_L\times [0,+\infty[/\cV_L\times\{0\}$ be the
cone over $\cV_L$; it is a conical vector bundle over $cL$. \\ 
Notice that the
cone  $cS^{N-1}$ over $S^{N-1}$ is isomorphic to $\RR^{N}_{\bullet}$ which is
$\RR^{N}$ with $0$ as a singular 
point. We will say that $cL$ is {\it embedded} in $\RR^{N}_{\bullet}$
and that $c\cV_L$ is the {\it tubular neighborhood of this embedding}.

\smallskip Now, let $X =(X, \fS, \cC)$ be a compact conical
stratification. Let 
$\cC=\{(\cN_s,\rho_s,L_s),\ s\in \fS \}$ be the set of control data, where $\cN_s$ is a cone over
$L_s$  and choose a trivialisation map $\varphi_s:\cN_s {\rightarrow}
cL_s$ for each singular
point $s$ . For $N\in
\NN$ large enough, one can find an embedding $j: \smx =
X\setminus\fS\rightarrow \RR^N$ such that : 
\begin{itemize} \item for any $s\in \fS$, $j\circ
\varphi_c^{-1}(L_s\times \{\lambda\})$ lies on a sphere
$S(O_s,\lambda)$ centered on $O_s$ and of radius $\lambda$ for
$\lambda \in ]0,1[$,
\item the open balls $B(O_s,1)$ centered on $O_s$ and of radius $1$
  are disjoint and,
\item for each singular
point $s$ there is an embedding $j_{L_s} : L_s \to S(O_c,1)\subset
\RR^N$ such that : \begin{equation*}
  \psi_s \circ  j\circ \varphi_c^{-1}\vert_{L\times ]0,1[ }
    = j_{L_s}\times \mbox{Id},
\end{equation*}
where $\psi_s : B(O_s,1)\setminus \{O_s\}
\rightarrow S(O_s,1)\times ]0,1[$ is the canonical diffeomorphism.
\end{itemize}
Let $\cV_{L_s} \to L_s$ be the 
normal bundle of the embedding $j_{L_s}$ and $\cV\rightarrow \smx$ be
the normal bundle of the embedding $j$. Then we can identify the
restriction of $\cV$ to $\cN_s\vert_{]0,1[}:=\{z\in \cN_s \vert 0<\rho_s(z)<1 \}$
  with $\cV_{L_s} \times ] 0, 1[$. Let $c\cV_{L_s}=\cV_{L_s} \times [
  0, 1[/_{\cV_{L_s} \times \{0\}}$ be the cone over $L_s$. We define
    the conical manifold  $$\cW =\cV \cup_{s\in \fS} 
c\cV_{L_s}$$ by glueing with $T\varphi_s$ the restriction of 
$\cV$ over $\cN_s\vert_{]0,1[}$ with $c\cV_{L_s}\setminus \{s\}$.
The conical manifold $\cW$ is a conical vector bundle over $X$. It
follows that $\cW$ is a sub-stratified pseudomanifold of
$({\RR}^N)^{\fS}$ which is $\RR^N$ with $\{O_s\}_{s\in \fS }$ as singular
  points. We will say
that $\cW$ is the {\it tubular neighborhood of the embedding of
$X$  in $(\RR^N)^{\fS}$}.

\section{Lie groupoids and their Lie algebroids\label{sec2}}

We refer to \cite{Re,CW,Mc} for the classical definitions and
construction related to groupoids and their Lie algebroids. 

\subsection{Lie groupoids} 

Groupoids, and especially differentiable groupoids will play an
important role in what follows, so we recall the basic
definitions and results needed for this
paper. Recall first that a {\em groupoid} is a small
category in which every morphism is an isomorphism.

Let us make the notion of a groupoid more explicit. Thus, a
groupoid $\cG$ is a pair $(\Gr0,\Gr1)$ of sets together  with
structural morphisms $u : \Gr0 \to \Gr1$, $s,r : \Gr1 \to \Gr0$,
$\iota : \Gr1 \to \Gr1$, and, especially, the multiplication $\mu$
which is defined for pairs $(g, h) \in \Gr1 \times \Gr1$ such that
$s(g) = r(h)$. Here, the set $\Gr0$ denotes the set of objects (or
units) of the groupoid, whereas the set $\Gr1$ denotes the set of
morphisms of $\cG$. Each object of $\cG$ can be identified with a
morphism of $\cG$, the identity morphism of that object, which
leads to an injective map $u : \Gr0 \to \cG$. Each morphism $g \in
\cG$ has a ``source'' and a ``range.'' We shall denote by $s(g)$
the {\em source} of $g$ and by $r(g)$ the {\em range} of $g$. The
inverse of a morphism $g$ is denoted by $g^{-1}=\iota(g)$. The
structural maps satisfy the following properties:
\begin{enumerate}[(i)]
\item  $r(gh)=r(g)$ and $s(gh)=s(h)$,
       for any pair $g,h$ satisfying $s(g) = r(h)$,
\item  $s(u(x))=r(u(x))=x$,  $u(r(g))g=g$,
$gu(s(g))=g$,
\item  $r(g^{-1})=s(g)$,\ $s(g^{-1})=r(g)$,\ $gg^{-1}=u(r(g))$,\
and $g^{-1}g=u(s(g))$,
\item  the partially defined multiplication $\mu$ is associative.
\end{enumerate}

We shall need groupoids with smooth structures.

\begin{definition}\label{Almost.differentiable}
A {\em Lie groupoid} is a groupoid
\begin{equation*}
        \cG=(\Gr0,\Gr1,s,r,\mu,u,\iota)
\end{equation*}
such that $\Gr 0$ and $\Gr1$ are manifolds with corners, the
structural maps $s,r,\mu,u,$ and $\iota$ are differentiable, the
domain map $s$ is a submersion  and $\cG_x :=
s^{-1}(x)$, $x\in M$, are all Hausdorff manifolds without corners.
\end{definition}

\smallskip The term ``differentiable groupoid'' was used in the past instead
of ``Lie groupoid,'' whereas ``Lie groupoid'' had a more
restricted meaning \cite{Mc}. The usage has changed however more
recently, and our definition reflects this change.

\smallskip An example of a Lie groupoid that will be used repeatedly below is
that of {\em pair groupoid}, which we now define. Let $M$ be a
smooth manifold. We let $\Gr0 = M$, $\Gr1 = M \times M$, $s(x, y)
= y$, $r(x, y) = x$, $(x, y)(y, z) = (x, z)$, and embedding $u(x)
= (x, x)$. The inverse is $\iota (x, y) = (y, x)$.

\smallskip The infinitesimal object associated to a Lie groupoid is its ``Lie
algebroid,'' which we define next.

\begin{definition}\label{Lie.Algebroid}
A  {\em Lie algebroid} $A$ over a manifold $M$ is a vector bundle
$A \to M$, together with a Lie algebra structure on the space
$\Gamma(A)$ of smooth sections of $A$ and  a bundle map $\varrho:
A \rightarrow TM$ whose extension to sections of these bundles
satisfies

(i) $\varrho([X,Y])=[\varrho(X),\varrho(Y)]$, and

(ii) $[X, fY] = f[X,Y] + (\varrho(X) f)Y$,

\noindent for any smooth sections $X$ and $Y$ of $A$ and any
smooth function $f$ on $M$.
\end{definition}

\smallskip The map $\varrho$ is called the {\em anchor map} of $A$. Note that
we allow the base $M$ in the definition above to be a manifold
with corners.

\smallskip The Lie algebroid associated to a differentiable groupoid $\cG$ is
defined as follows \cite{Mc}. The vertical tangent bundle (along
the fibers of $s$) of a differentiable groupoid $\cG$ is, as
usual,
\begin{equation}
        T_{vert} \cG = \ker s_*
        = \bigcup_{x\in M} T \cG_x \subset T\cG.
\end{equation}
Then $A(\cG) := T_{vert} \cG\big |_{M}$, the restriction of the
$s$-vertical tangent bundle to the set of units, defines the
vector bundle structure on $A(\cG)$.

\smallskip We now construct the bracket defining the Lie algebra structure on
$\Gamma(A(\cG))$. The right translation by an arrow $g \in \cG$
defines a diffeomorphism
\begin{equation*}
        R_g:\cG_{r(g)}\ni g' \longmapsto g'g \in \cG_{d(g)}.
\end{equation*}
A vector field $X$ on $\cG$ is called {\em $s$-vertical} if
$s_*(X(g)) = 0$ for all $g$. The $s$-vertical vector fields are
precisely the vector fields on $\cG$ that can be restricted to
vector fields on the submanifolds $\cG_x$. It makes sense then to
consider right--invariant vector fields on $\cG$. It is not
difficult to see that the sections of $A(\cG)$ are in one-to-one
correspondence with $s$--vertical, right--invariant vector fields
on $\cG$.

\smallskip The Lie bracket $[X,Y]$ of two $s$--vertical, right--invariant
vector fields $X$ and $Y$ is also $s$--vertical and
right--invariant, and hence the Lie bracket induces a Lie algebra
structure on the sections of $A(\cG)$. To define the action of the
sections of $A(\cG)$ on functions on $M$, let us observe that the
right invariance property makes sense also for functions on $\cG$,
and that $\cC^\infty(M)$ may be identified with the subspace of
smooth, right--invariant functions on $\cG$. If $X$ is a
right--invariant vector field on $\cG$ and $f$ is a
right--invariant function on $\cG$, then $X(f)$ will still be a
right invariant function. This identifies the action of
$\Gamma(A(\cG))$ on $\cC^\infty(M)$.

\subsection{Pull back groupoids}
Let $G\rightrightarrows M$ be a groupoid with source $s$ and range
$r$. If $f:N\rightarrow M$ is a surjective map, the {\it pull
back} groupoid $\pb{f}(G)\rightrightarrows N$ of $G$ by $f$ is by
definition the set
\begin{equation*}
\pb{f}(G):=\{(x,\gamma,y)\in N\times G\times N \ \vert \
r(\gamma)=f(x),\ s(\gamma)=f(y)\}
\end{equation*}
with the structural morphisms given by
\begin{enumerate}
\item the unit map $x \mapsto (x,f(x),x)$,

\item the source map $(x,\gamma,y)\mapsto y$ and range map
$(x,\gamma,y) \mapsto x$,

\item the product $(x,\gamma,y)(y,\eta,z)=(x,\gamma \eta ,z)$ and
inverse $(x,\gamma,y)^{-1}=(y,\gamma^{-1},x)$.
\end{enumerate}

\smallskip The results of  \cite{MRW} apply to show that the groupoids $G$
and $\pb{f}(G)$ are Morita equivalent.

\smallskip Let us assume for the rest of this subsection that $G$ is a smooth
groupoid and that $f$ is a surjective submersion, then $\pb{f}(G)$
is also a Lie groupoid. Let $(\cA(G), q, [\ ,\ ])$ be the Lie
algebroid of $G$ (which is defined since $G$ is smooth). Recall
that $q: \cA(G) \to TM$ is the anchor map. Let
$(\cA(\pb{f}(G)),p,[\ ,\ ])$ be the Lie algebroid of $\pb{f}(G)$
and $Tf : TN \to TM$ be the differential of $f$. Then we claim
that there exists an isomorphism
\begin{equation*}
    \cA(\pb{f}(G)) \simeq \{(V, U)\in TN \times \cA(G) \
    \vert \ Tf(V)=q(U) \in TM \}
\end{equation*}
under which the anchor map $p:\cA(\pb{f}(G)) \rightarrow TN$
identifies with the projection $TN \times \cA(G) \to TN$. In
particular, if $(U, V) \in \cA(\pb{f}(G))$ with $U \in T_xN$ and
$V \in \cA_y(G)$, then $y = f(x)$.

\subsection{Quasi-graphoid and almost injective Lie algebroid}

Our Lie groupoids arise mostly as Lie groupoids with a given
Lie algebroid. This is because often in Analysis, one is given the
set of derivations (differential operators), which forms a Lie
algebra under the commutator. The groupoids are then used to
``quantize'' the given Lie algebra of vector fields to algebra of
pseudodifferential operators \cite{aln2, MelroseScattering,
Monthubert, NWX}. This has motivated several works on the
integration of Lie algebroids \cite{CrainicFernandes, De,
NistorINT}. We recall here some useful results of the first named
author \cite{De} on the integration of some Lie algebroids. See
also \cite{CrainicFernandes, Mc, NistorINT}.

\begin{proposition}\label{propdef}
Let $G \overset{s}{\underset{r}{\rightrightarrows}} M$ be a Lie
groupoid over the manifold $M$. Let us denote by $s$ its domain
map, by $r$ its range map, and by $u:M \longrightarrow G$ its unit
map. The two following assertions are equivalent:

\noindent 1.\ If $\nu :V \longrightarrow G$ is a local section of
$s$ then $r \circ \nu = 1_V$ if, and only if, $\nu = u {\vert}_V$.

\noindent 2.\ If $N$ is a manifold, $f$ and $g$ are two smooth
maps from $N$ to $G$ such that:
\begin{enumerate}
\item[(i)] $s \circ f = s \circ g$ and $r \circ f= r \circ g$,

\item[(ii)] one of the following maps $s \circ f$ and $r \circ f$
is a submersion,
\end{enumerate}
then $f=g$.
\end{proposition}

\begin{definition} A Lie groupoid that satisfies
one of the two equivalent properties of Proposition \ref{propdef}
will be called a  {\it quasi-graphoid}.
\end{definition}

\smallskip Suppose that $G {\rightrightarrows} M$ is a quasi-graphoid and
denote by $\cA G=(p:\cA G \rightarrow TM,[\ ,\ ]_{\cA})$ its Lie
algebroid. A direct consequence of the previous definition is that
the anchor $p$ of $\cA G$ is injective when restricted to a dense
open subset of the base space $M$. In other words the anchor $p$
induces an injective morphism $\tilde{p}$ from the set of smooth
local sections of $\cA G$ onto the set of smooth local tangent
vector fields over $M$. In this situation we say that the Lie
algebroid $\cA G$ is {\it almost injective}.

\smallskip A less obvious remarkable property of a
quasi-graphoid is that its $s$-connected component is determined
by its infinitesimal structure. Precisely:

\begin{proposition}\cite{De}
Two $s$-connected quasi-graphoids having the same space of units
are isomorphic if, and only if, their Lie algebroids are
isomorphic.
\end{proposition}

Note that we are not requiring the groupoids in the above
proposition to be $s$-simply connected. The main result of
\cite{De} is the following:

\begin{theorem}\label{Integration} Every almost injective Lie algebroid
is integrable by an $s$-connected quasi-graphoid (uniquely by the
above proposition).
\end{theorem}

\smallskip Finally, let $\cA$ be a smooth vector bundle over a manifold $M$
and $p:\cA \rightarrow TM$ a morphism. We denote by $\tilde{p}$
the map induced by $p$ from the set of smooth local section of
$\cA$ to the set of smooth local vector fields on $M$. Notice that
if $\tilde{p}$ is injective then $\cA$ can be equipped with a Lie
algebroid structure over $M$ with anchor $p$ if, and only if, the
image of $\tilde{p}$ is stable under the Lie bracket.
\smallskip

\begin{examples} {\bf Regular foliation:}
A smooth regular foliation $\cF$ on a manifold $M$ determines an
integrable subbundle $F$ of $TM$. Such a subbundle is an (almost)
injective Lie algebroid over $M$. The holonomy groupoid of $\cF$
is the $s$-connected quasi-graphoid which integrates $F$
\cite{Wi}.\smallskip

\noindent{\bf Tangent groupoid:} \label{Tangent-groupoid} One
typical example of a quasi-graphoid is the tangent groupoid of A.
Connes \cite{CoBook}. Let us denote by $A \sqcup B$ the {\em
disjoint union} of the sets $A$ and $B$. If $M$ is a smooth
manifold, the tangent groupoid of $M$ is the disjoint union
\begin{equation*}
    \cG^t_M=TM\times \{0\} \sqcup M\times M \times ]0,1]
    \rightrightarrows M\times [0,1].
\end{equation*}
In order to equip $\cG^t_M$ with a smooth structure, we choose a
riemannian metric on $M$ and we require that the map
\begin{equation*}
\begin{array}{ccc}
    V\subset TM \times [0,1] & \longrightarrow & \cG_M^t
    \\ (x,V,t) & \mapsto &
    \left\{\begin{array}{ll} (x,V,0) \mbox{ if } t=0 \\
    (x,exp_x(-tV),t) \mbox{ if } t\not= 0 \end{array} \right.
\end{array}
\end{equation*}
be a smooth diffeomorphism onto its image, where $V$ is open in
$TM \times [0,1]$ and contains $TM\times \{0\}$. The tangent
groupoid of $M$ is the $s$-connected quasi-graphoid which
integrates the almost injective Lie algebroid:
\begin{equation*}\begin{array}{cccc} p^t_M : & \cA \cG_M^t = TM\times [0,1] &
\longrightarrow & T(M\times [0,1])\simeq TM \times T[0,1] \\ &
(x,V,t) & \mapsto & (x,tV;t,0)
\end{array}\end{equation*}
\end{examples}

\subsection{Deformation of quasi-graphoids} In this paper, we
will encounter deformation groupoids. The previous results give
easy aguments to get sure that these deformation groupoids can be
equipped with a smooth structure. For example, let $G_i
\rightrightarrows M$, $i=1,2$, be two $s$--connected
quasi--graphoids over the manifold $M$ and let $\cA G_i=(p_i:\cA
G_i \rightarrow TM , [\ , \ ]_{\cA_i})$ be the corresponding Lie
algebroid. Suppose that:
\begin{itemize}
  \item The bundles $\cA G_1$ and $\cA G_2$ are
  isomorphic,

\item There is a morphism $p:\cA:=\cA G_1 \times [0,1] \rightarrow
TM \times T([0,1])$ of the form:

\begin{equation*}
    p(V,0)=(p_1(V);0,0) \ \mbox{ and } \
    p(V,t)=(p_2 \circ \Phi(V,t);t,0) \mbox{ if } t\not=0,
\end{equation*}
where $\Phi:\cA G_1 \times ]0,1]\rightarrow \cA G_2 \times ]0,1]$
is an isomorphism of bundles over $M\times ]0,1]$. Moreover the
image of $\tilde{p}$ is stable under the Lie bracket.
\end{itemize}
In this situation, $\cA$ is an almost injective Lie algebroid that
can be integrated by the groupoid $H=G_1\times \{0\} \cup
G_2\times ]0,1]\rightrightarrows M\times [0,1]$. In particular,
there is a smooth structure on $H$ compatible with the smooth
structure on $G_1$ and $G_2$.

%%%%%%%%%%%%%%%%%%%%%%%%%%%%%%%%%
%%%%%%%%%%%%%%%%%%%%%%%%%%%%%%%%%%
%%%%%%%%%%%%%%%%%%%%%

\section{A non-commutative tangent space for conical
pseudomanifolds\label{sec3}}

In order to obtain an Atiyah-Singer type topological index theorem
for our conical pseudomanifold $X$, we introduce in this chapter a
suitable notion of tangent space to $X$ and a suitable normal
space to an embedding of $X$ in $\RR^{N+1}$ that sends the
singular point to $0$ and $\smx$ to $\{x_1 > 0\}$.

\subsection{The $\fS$-tangent space and the tangent groupoid of a
conical space} We recall here a construction from \cite{DL} that
associates to a conical pseudomanifold $X$ a groupoid $T^\fS X$
that is a replacement of the notion of tangent space of $X$ (for the
purpose of studying $K$-theory) in the sense the $C^*$-algebras
$C^*(T^\fS X)$ and $C(X)$ are $K$-dual \cite{DL}.

\smallskip \noindent Let $(X,\fS,\cC)$ be a conical
pseudomanifold. Without loss of generality, we can assume that $X$
has only one singular point. Thus $\fS=\{c\}$ is a single point
and $\cC=\{(\cN,\rho,L)\}$, where $\cN \simeq cL$ is a cone over
$L$ and $\rho$ is the defining function of the cone. We set
$\rho=+\infty$ outside $\cN$. We let $\smx=X\setminus \{c\}$.
Recall that $\smx$ is a smooth manifold. We denote by $O_X$ the
open set $O_X=\{ z\in \smx \ \vert \ \rho(z)< 1 \}$.

\smallskip \noindent
At the level of sets, the {\it $\fS$-tangent space} of $X$ is the
groupoid: \begin{equation*}T^{\fS} X:= T\smx \vert_{\smx\setminus
O_X} \sqcup O_X\times O_X \rightrightarrows \smx \
.\end{equation*}

\noindent Here, the groupoid $T\smx \vert_{\smx\setminus
O_X}\rightrightarrows \smx\setminus O_X $ is the usual tangent
vector bundle $T\smx$ of $\smx$ restricted to the closed subset
$\smx\setminus O_X=\{z\in \smx \ \vert \ \rho(z)\geq 1 \}$. The
groupoid $O_X\times O_X \rightrightarrows O_X$ is the pair
groupoid over $O_X$.

\smallskip \noindent
The {\it tangent groupoid} of $X$ is, as in the regular case
\cite{CoBook}, a deformation of its ``tangent space'' to the pair
groupoid over its units:
\begin{equation*}
    \cG^t_X:= T^{\fS}X \times \{0\} \sqcup
    \smx \times \smx \times
    ]0,1] \rightrightarrows \smx\times [0,1].
\end{equation*}
Here, the groupoid $\smx \times \smx \times ]0,1]
\rightrightarrows \smx\times ]0,1]$ is the product of the pair
groupoid on $\smx$ with the set $]0,1]$.

\smallskip In order to equip $\cG^t_X$, and so $T^{\fS}X$,
with a smooth structure we have to choose a {\it glueing
function}. First choose a positive smooth map $\tau: \RR
\rightarrow \RR$ such that $\tau([0,+\infty[)=[0,1]$,
$\tau^{-1}(0)=[1,+\infty[$ and $\tau'(t)\not=0$ for $t < 1$. We
denote by $\tau_X: X\rightarrow \RR$ the map which assigns $\tau
(\rho(x))$ to $x\in \smx\cap \cN$ and $0$ elsewhere. Thus
$\tau_X(\smx)=[0,1 [$, $\tau_X$ restricted to $O_X =\{z\in \cN \
\vert \ 0<\rho(z) < 1\}$ is a submersion and $\tau_X^{-1}(0)=\smx
\setminus O_X$.

\begin{proposition} \cite{DL}
There is a unique structure of Lie groupoid on $\cG_X^t$ such
that its Lie algebroid is the bundle $T\smx\times [0,1]$ whith
anchor $p:(x,V,t)\in T\smx\times [0,1] \mapsto
(x,(t+\tau_X^2(x))V;t,0) \in T\smx\times T[0,1]$.
\end{proposition}

\smallskip Let us notice that the map $p$ is injective when restricted to
$\smx\times ]0,1]$, which is a dense open subset of $\smx\times
[0,1]$. Thus there exists one, and only one, structure of (almost
injective) Lie algebroid on $T\smx\times [0,1]$ with $p$ as anchor
since the family of local vector fields on $\smx$ induced by the
image by $p$ of local sections of $T\smx\times [0,1]$ is stable
under the Lie bracket. We know from \cite{De, NistorINT} that such
a Lie algebroid is integrable.  Moreover, according to theorem
\ref{Integration}, there is a unique Lie groupoid  which
integrates this algebroid and restricts over $\smx\times ]0,1]$ to
$\smx\times \smx \times ]0,1] \rightrightarrows \smx\times ]0,1]$.

\medskip \noindent  Let us give an
alternative proof of the previous proposition.

\noindent \begin{proof} Recall that the (classical) tangent
groupoid of $\smx$ is
\begin{equation*}
    \cG^t_{\smx}=T\smx\times \{0\} \sqcup
    \smx\times \smx \times ]0,1]
    \rightrightarrows \smx\times [0,1]
\end{equation*}
and that its Lie algebroid is the bundle $T\smx\times [0,1]$ over
$\smx\times [0,1]$ with anchor $(x,V,t)\in T\smx\times [0,1]
\mapsto (x,tV,t,0) \in T\smx\times T[0,1]$. Similary, one can
equip the groupoid $H=T\smx\times \{(0,0)\} \sqcup \smx \times
\smx \times [0,1]^2\setminus \{(0,0)\}$ with a unique smooth
structure such that its Lie algebroid is the bundle $T\smx\times
[0,1]^2$ with anchor the map \begin{equation*}\begin{array}{cccc}
  p: & \cA=T\cN_1\times [0,1]\times [0,1] & \rightarrow &
  T\cN_1\times
  T([0,1]) \times T([0,1]) \\ & (x,V,t,l) & \mapsto &
  (x,(t+l)V;t,0;l,0)
\end{array}\end{equation*}

\smallskip Let $\delta:H \rightarrow \RR$ be the map which sends any $\gamma
\in H$ with source $s(\gamma)=(y,t,l)$ and range
$r(\gamma)=(x,t,l)$ to $\delta(\gamma)=l-\tau_X(x)\tau_X(y)$. One
can check that $\delta$ is a smooth submersion, so $H_{\delta} :=
\delta^{-1}(0)$ is a submanifold of $H$. Moreover 
$H_{\delta} := \delta^{-1}(0)$ inherits from $H$ a structure of
Lie groupoid over $\smx\times [0,1]$ whose Lie algebroid is
given by
\begin{equation*}
\begin{array}{ccc} T\smx \times [0,1] & \rightarrow & T\smx
\times
  T([0,1]) \\ (x,V,t) & \mapsto &
  (x,(t+\tau_X^2(x))V;t,0) \end{array}
\end{equation*}
The groupoid $H_{\delta}$ is (obviously isomorphic) to $\cG^t_X$.
\end{proof}

\smallskip We now introduce the tangent groupoid of a stratified
pseudomanifold.

\begin{definition}\label{def.tangent.groupoid}
The groupoid $\cG^t_X$ equipped with the smooth structure associated
with a glueing function $\tau$ as above is called a {\em tangent
groupoid} of the stratified pseudomanifold $(X,\fS,\cC)$. The
corresponding $\fS$-tangent space is the groupoid
$T^{\fS}X\simeq\cG^t_X \vert_{\smx\times \{0\}}$ equipped with the
induced smooth structure.
\end{definition}

\begin{remarks}\label{miscellaneous-properties-TSX}
We will need the following remarks. See \cite{DL} for a proof.
\begin{enumerate}[(i)]
\item If $X$ has more than one singular point, we let, for any $s\in
\fS$,
\begin{equation*}
    O_s:=\{z\in \smx\cap \cN_s \ \vert \ \rho_s(z)<1 \},
\end{equation*}
and we define $O=\sqcup_{s\in \fS}O_s$. The {\em $\fS$-tangent
space to $X$} is then
\begin{equation*}
    T^{\fS}X:= T\smx\vert_{\smx\setminus O} \sqcup_{s\in \fS}
    O_s\times O_s \rightrightarrows \smx,
\end{equation*}
with the analogous smooth structure. In this situation the Lie
algebroid of $\cG^t_X$ is defined as previously with
$\tau_X:X\rightarrow \RR$ being the map which assigns
$\tau(\rho_s(z))$ to $z\in \smx \cap \cN_s$ and $0$ elsewhere.
\item The orbit space of $T^{\fS}X$ is
topologicaly equivalent to $X$: there is a canonical isomorphism
between the algebras $C(X)$ and $C(X/T^{\fS}X)$.
\item The tangent groupoid and the $\fS$-tangent space depend on
the glueing. Nevertheless the $K$-theory of the $C^*$-algebras
$C^*(\cG^t_X)$ and $C^*(T^{\fS}X)$ do not.
\item The groupoid $T^\fS X$ is a continuous field of amenable
groupoids parametrized by $X$, thus $T^\fS X$ is amenable as well.
It follows that $\cG^t_X$ is also amenable as a continuous field
of amenable groupoids parametrised by $[0,1]$. Hence the reduced
and maximal $C^*$-algebras of $T^\fS X$ and of $\cG^t_X$ are equal
and they are nuclear.
\end{enumerate}
\end{remarks}

\begin{examples}
Here are two basic examples.
\begin{enumerate}[(i)]
\item When $X$ is a smooth manifold, that is $X_0=\emptyset$ and
$\smx=X$, the previous construction gives rise to the usual
tangent groupoid
\begin{equation*}
    \cG^t_X=TX\times \{0\} \sqcup X\times X \times ]0,1]
    \rightrightarrows X \times [0,1].
\end{equation*}
Moreover, $T^{\fS}X = TX \rightrightarrows X$ is the usual tangent
space.\
\item Let $L$ be a manifold and consider the (trivial) cone
$cL=L\times[0,+\infty [/L\times\{0\}$ over $L$. In this situation
$\smx=L \times ]0,+\infty[$, $O_X=L\times ]0,1[$ and
\begin{equation*}
    T^{\fS} X =\, T\big(L \times [1,+\infty[\,\big) \sqcup\,
    \underbrace{\,L\times\, ]0,1[\, \times\,
    L\, \times\, ]0,1[\,}_{\mbox{the pair groupoid}}\
    \rightrightarrows\ L\, \times\, ]0,+\infty[ \ ,
\end{equation*}
where $T (L\times [1,+\infty[)$ denotes the restriction to
$L\times [1,+\infty[$ of the tangent space $T(L\times \RR)$. The
general case is always locally of this form.
\end{enumerate}
\end{examples}

\subsection{The deformation groupoid of a conical vector
bundle\label{deformation}} Let $(E,\fS_E,\cC_E)$ be a conical
vector bundle over $(X,\fS_X,\cC_X)$ and denote by $\pi:
E\rightarrow X$ the corresponding projection. From the definition,
$\pi$  restricts to a smooth vector bundle map $\pi^{\circ} :
\smE \rightarrow \smx$. We let $\pi_{[0,1]}=\pi^{\circ} \times id
: \smE \times [0,1] \rightarrow \smx \times [0,1]$.

\smallskip We consider the tangent groupoids
$\cG^t_X \rightrightarrows \smx$ for $X$ and $\cG^t_E
\rightrightarrows \smE$ for $E$ equipped with a smooth structure
constructed using the same glueing fonction $\tau$ (in particular
$\tau_X \circ \pi = \tau_E$). We denote by
$\pb{\pi_{[0,1]}}(\cG^t_X)\rightrightarrows \smE \times [0,1]$ the
pull back of $\cG^t_X$ by $\pi_{[0,1]}$.

\smallskip  
Our next goal is to associate to the conical vector bundle
$E$ a deformation groupoid $\cT^t_{E}$ using
$\pb{\pi_{[0,1]}}(\cG^t_X)$ to $\cG^t_E$.  More precisely, we
define:
\begin{equation*}
    \cT^t_{E}:= \cG^t_E \times \{0\} \sqcup
    \pb{\pi_{[0,1]}}(\cG^t_X)\times ]0,1]
    \rightrightarrows \smE \times [0,1] \times [0,1].
\end{equation*}
In order to equip $\cT^t_E$ with a smooth structure, we first
choose a smooth projection $P:T \smE \rightarrow
\mbox{Ker}(T\pi)$.

\smallskip A simple calculation shows that the Lie algebroid of
$\pb{\pi_{[0,1]}}(\cG^t_X)$ is isomorphic to the bundle $T
\smE \times [0,1]$ endowed with the almost injective anchor map
\begin{equation*}
    (x,V,t) \mapsto (x,P(x,V)+(t+\tau_E(x)^2)(V-P(x,V));t,0).
\end{equation*}
We consider the bundle $\cA=T \smE \times [0,1] \times [0,1]$ over
$\smE\times [0,1] \times [0,1]$ and the almost injective morphism:
\begin{equation*}
\begin{array}{cccc}
    p: & \cA =T \smE \times [0,1] \times [0,1]
    & \rightarrow & T \smx \times T[0,1] \times T[0,1] \\
    & (x,V,t,l) & \mapsto & (x,(t+\tau_E^2(x))V+lP(x,V)).
\end{array}
\end{equation*}
The image of $\tilde{p}$ is stable under the Lie bracket, thus
$\cA$ is an almost injective Lie algebroid. Moreover, the
restriction of $\cA$ to $\smE \times [0,1] \times \{0\}$ is the
Lie algebroid of $\cG_E^t$ and its restriction to $\smE \times
[0,1] \times ]0,1]$ is isomorphic to the Lie algebroid of
$\pb{\pi_{[0,1]}}(\cG^t_X) \times ]0,1]$. Thus $\cA$ can be
integrated by $\cT^t_E$. In particular, $\cT^t_E$ is a smooth
groupoid. In conclusion, the restriction of $\cT^t_E$ to
$\smE \times \{0\} \times [0,1]$ leads to a Lie groupoid:
\begin{equation*}
    \cH_{E} = T^{\fS}E\times \{0\} \sqcup
    \pb{\pi}(T^{\fS}X) \times ]0,1]
    \rightrightarrows \smE \times [0,1],
\end{equation*}
called a {\it Thom groupoid} associated to the conical vector bundle
$E$ over $X$.

\smallskip The following example explains what these constructions
become if there are no singularities.

\begin{example}
Suppose that $p:E\rightarrow M$ is a smooth vector bundle over the
smooth manifold $M$. Then $\fS_E = \fS_M=\emptyset$, $\cG^t_E =
TE\times \{0\} \sqcup E \times E \times ]0,1] \rightrightarrows E
\times [0,1]$ and $\cG^t_M = TM \times \{0\} \sqcup M \times M
\times ]0,1] \rightrightarrows M \times [0,1]$ are the usual
tangent groupoids. In these examples associated to a smooth
vector bundle, $\tau_E$ is the zero map. The groupoid $\cT^t_{E}$ will
then be given by
\begin{equation*}
    \cT^t_{E}= TE\times \{0 \}\times \{0 \}
    \sqcup \pb{p}(TM)\times \{0\} \times ]0,1]
    \sqcup E\times E \times ]0,1] \times [0,1]
    \rightrightarrows E\times [0,1]\times [0,1]
\end{equation*}
and is smooth. Similarly, the Thom groupoid will be given
by: $\cH_E:= TE\times \{0 \} \sqcup \pb{p}(TM) \times
]0,1]\rightrightarrows E\times [0,1]$.
\end{example}

\smallskip We now return to the general case of a conical vector
bundle.

\begin{remark}
The groupoids $\cT_E$ and  $\cH_E$ are
continuous fields of amenable groupoids parametrized by $[0,1]$.
Thus they are amenable, their reduced and maximal $C^*$-algebras
are equal, and are nuclear.
\end{remark}

%%%%%%%%%%%%%%%%%%%%%%%%%%%%%%%%%%
%%%%%%%%%%%%%%%%%%%%%%%%%%%%%%%%%%
%%%%%%%%%%%%%%%%%%%%

\section{The analytical index\label{sec4}}

Let $X$ be a conical pseudomanifold, and let
\begin{equation*}
    \cG_X^t=\smx\times \smx \times ]0,1]
    \sqcup T^{\fS} X \times \{0\}
    \rightrightarrows \smx\times [0, 1]
\end{equation*}
be the tangent groupoid (unique up to isomorphism) for $X$ for a given
glueing function. Also, let $T^\fS X \rightrightarrows \smx$ be the
corresponding $\fS$-tangent space.

\smallskip Since the groupoid $\cG_X^t$ is a deformation groupoid of amenable
groupoids, it defines a $KK$-element \cite{DL, HS1}. More
precisely, let
\begin{equation*}
    e_1 : C^*(\cG_X^t) \rightarrow
    C^*(\cG^t_X \vert_{\smx\times\{1\}}) = \mathcal{K}(L^2(\smx))
\end{equation*}
be the evaluation at 1 and let $[e_1] \in KK(C^*(\cG_X^t),
\mathcal{K}(L^2(\smx))$ the element defined by $e_1$ in Kasparov's
bivariant $K$--theory. Similarly, the evaluation at 0 defines a
morphism $e_0: C^*(\cG^t_X) \rightarrow
C^*(\cG_X^t\vert_{\smx\times\{0\}}) = C^*(T^{\fS}X)$ and then an
element $[e_0] \in KK(C^*(\cG^t_X), C^*(T^{\fS}X))$. The kernel of
$e_0$ is contractible and so $e_0$ is $KK$-invertible. We let:
\begin{equation*}\widetilde{\partial} =[e_0]^{-1} \otimes [e_1] \in
KK(C^*(T^{\fS}X),\mathcal{K} ), \end{equation*} be the Kasparov
product over $C^*(\cG^t_X)$ of $[e_1]$ and the $K$-inverse of
$[e_0]$. Take $b$ to be a generator of
$KK(\mathcal{K},\mathbb{C})\simeq \mathbb{Z}$. We set $\partial =
\widetilde{\partial}\otimes b$. The element $\partial$ belongs to
$KK(C^*(T^{\fS}X),\CC)$.

\begin{definition} The map $(e_0)_* : K_0(C^*(\cG^t_X)) \to
K_0(C^*(T^{\fS}X))$ is an isomorphism and we define the {\em
analytical index map} by 
\begin{equation}
        \ind^X_a := (e_1)_* \circ (e_0)_*^{-1} :
        K_0(C^*(T^{\fS}X))\to K_0(\cK)\simeq\ZZ,
\end{equation}
or in other words, as the map defined by the Kasparov product with
$\partial$.
\end{definition}

\begin{remarks} {1.} Notice that in the case of a smooth manifold with
  the usual definition of tangent space and tangent groupoid, this
  definition leads to the classical definition of the analytical index
  map (\cite{CoBook} II.5).\\
{2.} One can associate to a Lie groupoid a different analytical
map. More precisely, when $G\rightrightarrows M$ is smooth, one
can consider the adiabatic groupoid which is a deformation
groupoid of $G$ on its Lie algebroid $\cA G$ \cite{NWX}:
\begin{equation*}
\cG^t:= \cA G \times \{0\} \sqcup G \times ]0,1]
\rightrightarrows M\times [0,1] \ .
\end{equation*}
Under some asumption $\cG^t$ defines a $KK$-element in $KK(C^*(\cA
G),C^*(G))$ and thus a map from $K_0(C^*(\cA G))$ to $K_0(C^*(G))$.
\end{remarks}

\smallskip Now, let $X$ be a conical pseudomanifold and $\fS$
  its set of singular points. Choose a singular point $s\in \fS$. Let
  us denote
$X_{s,+}:=\smx\setminus O_s$. The $\fS$-tangent space of $X$ is then
\begin{equation*}
    T^\fS X= O_s\times O_s \sqcup T^\fS X_{s,+} \rightrightarrows \smx
\end{equation*}
where $O_s\times O_s\rightrightarrows O_s$ is the pair groupoid
and $T^{\fS}X_{s,+}:= T^\fS X\vert_{X_{s,+}}$. Then we have the following exact
sequence of $C^*$-algebras:
\begin{eqnarray} \label{suite1} \begin{CD} 0\rightarrow
    \underbrace{C^*(O_s\times  O_s)}_{=\ \cK(L^2(O_s))} @>i>>
    C^*(T^{\fS}X) @>r_+>> C^*(T^{\fS} X_{s,+}) \rightarrow 0, \end{CD}
\end{eqnarray}
where $i$ is the inclusion morphism and $r_+$ comes from the
restriction of functions.
\begin{proposition} \label{shortsequence}
The exact sequence \eqref{suite1} induces the short exact sequence
\begin{equation*}
\begin{CD}
    0\rightarrow K_0(\cK)\simeq \ZZ @>i_*>>
    K_0(C^*(T^{\fS}X)) @>(r_+)_*>>
    K_0(C^*(T^{\fS }X_{s,+})) \rightarrow 0.
\end{CD}
\end{equation*}
Moreover  $\ind^X_a\circ i_*=\mbox{Id}_\ZZ$, thus
\begin{equation*}
    \big(\, \ind^X_a,\,(r_+)_*\ \big):K_0(C^*(T^{\fS}X)) \rightarrow
    \ZZ\oplus K_0(C^*(T^{\fS} X_{s,+}))
\end{equation*}
is an isomorphism .
\end{proposition}
\begin{proof}
In order to prove the first statement, let us first consider the
six terms exact sequence associated to the exact sequence of
$C^*$-algebra of \ref{suite1}. Then recall that $K_1(\cK)=0$. It
remains to show that $i_*$ is injective. This point is a
consequence of the second statement which is proved here'after.
\smallskip  Let $\cG_X^t \rightrightarrows \smx\times [0,1]$
be a tangent groupoid for $X$. Its restriction
$\cG_X^t\vert_{O_s\times[0,1]}$ to $O_s\times[0,1]$ is isomorphic
to the groupoid $(O_s\times O_s)\times[0,1]\rightrightarrows
O_s\times [0,1]$, the pair groupoid of $O_s$ parametrized by
$[0,1]$. The inclusion  of $C_0(\cG_X^t\vert_{O_s\times[0,1]})$ in
$C_0(\cG_X^t)$ induces a morphism of $C^*$-algebras
\begin{equation*}
    i^t:C^*(\cG_X^t\vert_{O_s\times[0,1]})
    \simeq \cK(L^2(O_s)) \otimes C([0,1])
    \rightarrow C^*(\cG_X^t).
\end{equation*}
Moreover, we have the following commutative diagram of
$C^*$-algebras morphisms:
\begin{equation*}
\begin{CD}
    C^*(T^\fS X) @<e_0<< C^*(\cG_X^t) @>e_1>> \cK(L^2(\smx))
    \\ @AiAA  @AAi^tA  @AAi_{\cK}A \\ \cK (L^2(O_s)) @<<ev_0<
    \cK(L^2(O_s)) \otimes C([0,1]) @>>ev_1>\cK (L^2(O_s))
\end{CD}
\end{equation*}
where $i_{\cK}$ is the isomorphism induced by the inclusion of the
pair groupoid of $O_s$ in the pair groupoid of $\smx$, and $ev_0$,
$ev_1$ are the evaluations map at $0$ and $1$. The $KK$-element
$[ev_0]$ is invertible and
\begin{equation*}
    [ev_0]^{-1} \otimes [ev_1] = 1
    \in KK(\cK(L^2(O_s)),\cK(L^2(O_s))).
\end{equation*}
Moreover $\cdot\otimes [i_\cK]$ induces an isomorphism from
$KK(\CC,\cK(L^2(O_s)))\simeq \ZZ$ onto
$KK(\CC,\cK(L^2(\smx)))\simeq \ZZ$.  Thus $[i]\otimes \tilde
\partial = [i_\cK]$, which proves that $\ind^X_a\circ
i_*=\mbox{Id}_\ZZ$ and ensures that $i_*$ is injective.
\end{proof}

%%%%%%%%%%%%%%%%%%%%%%%%%%%%%%%%
%%%%%%%%%%%%%%%%%%%%%%%%%%%%%%%%%
%%%%%%%%%%%%%%%%%%%%%%%

\section{The inverse Thom map\label{Thominv}}

Let $(E,\fS_E,\cC_E)$ be a conical vector bundle over
$(X,\fS_X,\cC_X)$ and $\pi: E\rightarrow X$ the corresponding
projection. We let
\begin{equation}\label{eq.elThom}
    \cH_{E}:= T^{\fS}E\times \{0\}
    \sqcup \pb{\pi}(T^{\fS}X) \times ]0,1]
    \rightrightarrows \smE \times [0,1]
\end{equation}
be the Thom groupoid of $E$, as before. The $C^*$-algebra of
$\cH_{E}$ is nuclear as well as the $C^*$-algebra of $T^{\fS}E$.
Thus $\cH_{E}$ defines a $KK$-element:
\begin{equation}\label{eq.invThom}
    \partial_{\cH_E} := [\epsilon_0]^{-1} \otimes [\epsilon_1]
    \in KK(C^*(T^\fS E), C^*(T^\fS X) ),
\end{equation}
where $\epsilon_1:C^*(\cH_E)\rightarrow
C^*(\cH_E\vert_{\smE\times\{1\}}) = C^*(\pb{\pi} (T^{\fS}X))$ is
the evaluation map at $1$ and $\epsilon_0: C^*(\cH_E) \rightarrow
C^*(\cH_E\vert_{\smE \times \{0\}}) = C^*(T^{\fS}E)$, the
evaluation map at $0$ is $K$-invertible.

\begin{definition} \label{def.invThom}
The element $\partial_{\cH_E} \in KK(C^*(T^\fS E), C^*(T^\fS X) )$
defined by Equation \eqref{eq.invThom} will be called the
{\it inverse Thom element}.
\end{definition}

\begin{dfprop} \label{Tinv}
Let $\cM$ be the isomorphism induced by the Morita equivalence
between $T^{\fS} X$ and $\pb{\pi}(T^{\fS} X)$ and let $\cdot
\otimes \partial_{\cH_E}$ be the right Kasparov product by
$\partial_{\cH_E}$ over $C^*(T^\fS E)$. Then the following diagram
is commutative:
\begin{equation*}
\begin{CD}
    K(C^*(T^{\fS} E)) @>{\ind_a^E}>> \ZZ \\
    @V{\cdot \otimes \partial_{\cH_E}}VV
    @AA{\ind_a^X}A \\ K(C^*(\pb{\pi}(T^{\fS} X)))
    @> \simeq>{\cM}> K(C^*(T^{\fS} X)).
\end{CD}
\end{equation*}
The map $T_{inv} := \cM \circ (\cdot \otimes
\partial_{\cH_E})$ is called {\em the inverse Thom map}.
\end{dfprop}

\begin{proof}
First consider the deformation groupoid $\cT^t_E$:
\begin{equation*}
    \cT^t_E:=\cG_E^t\times\{0\} \sqcup \pb{\pi_{[0,1]}}
    (\cG^t_X) \times ]0,1] \rightrightarrows
    \smE\times[0,1]\times[0,1] \ .
\end{equation*}
One can easily see that
\begin{equation*}
\begin{array}{lll}
    \cT^t_E & = & T^{\fS}E\times\{0\}\times \{0\}
    \sqcup \pb{\pi}(T^{\fS}X)\times \{0\}\times ]0,1]
    \sqcup \smE \times \smE \times]0,1] \times [0,1]\\
    & \simeq & \cH_{E} \times \{0\} \sqcup
    (\smE \times \smE \times[0,1])\times ]0,1] \, .
\end{array}
\end{equation*}

\smallskip The groupoid $\cT^t_E$ is equipped with a smooth structure
compatible with the smooth structures of $\cG_E^t\times\{0\}$,
$\pb{\pi_{[0,1]}}(\cG^t_X) \times ]0,1]$ as well as with the
smooth structures of $\cH_{E}$ and $\big(\smE \times \smE
\times[0,1]\big)\, \times\, ]0,1]$.

\smallskip 
We therefore have the following commutative diagram of evaluation
morphisms of $C^*$-algebras of groupoids:
\begin{equation*}
\begin{CD}
    C^*(\smE \times \smE ) @<{ev_0}<<
    C^*(\smE \times \smE \times [0,1])
    @>{ev_1}>> C^*(\smE \times \smE ) \\
    @A{e_1^E}AA @AA{q_{1,\cdot}}A
    @AA{\pb{\pi} e^X_1}A \\
    C^*(\cG^t_E)  @<{q_{\cdot,0}}<< C^*(\cT^t_E)
    @>{q_{\cdot,1}}>> C^*(\pb{\pi_{[0,1]}}(\cG^t_X)) \\
    @V{e_0^E}VV @VV{q_{0,\cdot}}V
    @VV{\pb{\pi}e^X_0}V \\
    C^*(T^\fS E) @<<{\epsilon_0}< C^*(\cH_E) @>>{\epsilon_1}>
    C^*(\pb{\pi}(T^\fS X))
\end{CD}
\end{equation*}

\noindent In this diagram, the $KK$-elements $[e_0^E]$,
$[\pb{\pi}e_0^X]$, $[q_{\cdot,0}]$, $[\epsilon_0]$, $[ev_1]$ and
$[ev_0]$ are invertible. Let $\cM:K(C^*(\pb{\pi}(T^\fS
X))\rightarrow K(C^*(T^\fS X))$ be the isomorphism induced by the
Morita equivalence between $\pb{\pi}(T^\fS X)$ and $T^\fS X$.
Also, let $x$ belong to $K(C^*(\pb{\pi}(T^\fS X))) =
KK(\CC,\pb{\pi}(T^\fS X))$. Then one can easily check the equality
\begin{equation*}
    \cM (x)\otimes \tilde{\partial}
    = x\otimes [\pb{\pi}e^X_0]^{-1}
    \otimes [\pb{\pi} e^X_1].
\end{equation*}
Of course $[ev_0]^{-1}\otimes [ev_1]=1 \in KK(C^*(\smE \times \smE
),C^*(\smE \times \smE )$. Thus the previous diagram implies that
for any $x\in K(C^*(T^\fS E))=KK(\CC,C^*(T^\fS E))$ we have:

\begin{equation*}\begin{array}{lll}
    \ind_a^X \circ T_{inv}(x) & = & x \otimes [\epsilon_0]^{-1}
    \otimes [\epsilon_1] \otimes [\pb{\pi}e^X_0]^{-1}
    \otimes [\pb{\pi} e^X_1]\otimes b \\
    & = & x \otimes [e^E_0]^{-1}\otimes [e^{E_1}]
    \otimes [ev_0]^{-1} \otimes[ev_1]\otimes b \\
    & = & \ind_a^E(x) \end{array}
\end{equation*}
\end{proof}

%%%%%%%%%%%%%%%%%%%%%%%%%%%%%%%%
%%%%%%%%%%%%%%%%%%%%%%%%%%%%%%%%
%%%%%%%%%%%%%%%%%%%%%%%%

\section{Index theorem\label{sec6}}

In this section, we state and prove our main theorem, namely, a
topological index theorem for conical pseudomanifolds in the setting
of groupoids. We begin with an account of the classical Atiyah-Singer
topological index theorem in our groupoid setting.

\subsection{A variant of the proof of Atiyah-Singer index
theorem for compact manifolds using groupoids} Let $\cV$ be the
normal bundle of an embedding of a smooth manifold $M$ in some
euclidean space. In this subsection, we shall first justify the
terminology of ``inverse Thom map'' we introduced for the map
$T_{inv}$ of Proposition \ref{Tinv} by showing that it coincides
with the inverse of the classical Thom isomorphism when $E=T\cV$
and $X=TM$.

\smallskip In fact, we will define the Thom
isomorphism when $X$ is a locally compact space and $E = N
\otimes \CC$ is the complexification of a real vector bundle $N
\to X$. As a consequence, we will derive a simple proof of the
Atiyah-Singer index theorem for closed smooth manifolds. Our
approach has the advantage that it extends to the singular
setting.

\smallskip  Let us recall some classical facts \cite{AS1,At}. If $p: E\to
X$ is a complex vector bundle over a locally compact space $X$,
one can define a Thom map
\begin{equation}\label{eq.Thom}
    i_!: K^0(X)\to K^0(E),
\end{equation}
which turns to be an isomorphism. This Thom map is defined as
follows. Let $x\in K^0(X)$ be represented by
$[\xi_0;\xi_1;\alpha]$ where $\xi_0,\xi_1$ are complex vector
bundles over $X$ and $\alpha : \xi_0\to\xi_1$ is an isomorphism
outside a compact subset of $X$. With no lost of generality, one
can assume that $\xi_0,\xi_1$ are hermitian and that $\alpha$ is
unitary outside a compact subset of $X$.

\smallskip  Let us consider next the endomorphism of the vector
bundle $p^*(\Lambda E)\to E$ given by
 \begin{equation*}
   (C\omega)(v) = C(v)\omega(v)=\frac{1}{\sqrt{1+\|v\|^2}}
   \left(v\wedge \omega(v) - v^*\llcorner \omega(v)\right)
 \end{equation*}
The endomorphism $C\omega$ is selfadjoint, of degree $1$ with
respect to the $\ZZ_2$--grading $\Lambda_0 =
\Lambda^{even}E,\Lambda_1=\Lambda^{odd}E$ of the space of exterior
forms. Moreover, we have that $(C\omega)^2 \to 1$ as $\omega$
approaches infinity in the fibers of $E$. Then, as we shall
see in the next Proposition, the Thom morphism $i_!$ of Equation
\eqref{eq.Thom} can be expressed, in terms of the Kasparov
products, as
\begin{equation*}
    i_!(x) := \left \lbrack \,
    \xi_0 \otimes \Lambda_0 \oplus \xi_1 \otimes \Lambda_1;\
    \xi_0 \otimes \Lambda_1 \oplus \xi_1 \otimes \Lambda_0;\
    \theta=
    \begin{pmatrix}
    N(1\otimes C) & M(\alpha^*\otimes 1) \\
    M(\alpha\otimes 1)  & -N(1 \otimes C)
    \end{pmatrix}\,
    \right\rbrack
\end{equation*}
where $M$ and $N$ are the multiplication operators by the
functions $M(v)=\frac{1}{\|v\|^2+1}$ and  $N=1-M$, respectively.

\begin{proposition} \label{propThom}
Let $p : E\To X$ be a complex vector bundle over a locally compact
base space $X$ and  $i_!: K^0(X)\To K^0(E)$ the corresponding Thom
map. Denote by $T$ the Kasparov element
\begin{equation*}
    T := \big(C_0(E,p^*(\Lambda E)),\,\rho\, ,\,
    C\,\big) \in KK(C_0(X),C_0(E))
\end{equation*}
where  $\rho$ is multiplication by functions. Then $i_!(x) = x \otimes T$ for any
$x\in K^0(X)$.
\end{proposition}

\begin{proof}
The isomorphism $K^0(X)\simeq KK(\CC,C_0(X))$ is such that to the
triple $[\xi_0;\xi_1;\alpha]$ there corresponds to the Kasparov
module:
\begin{equation*}
     x= \left( C_0(X,\xi),1,\widetilde{\alpha} \right), \quad
    \xi=\xi_0\oplus\xi_1\ \text{ and }\
    \widetilde{\alpha}=\begin{pmatrix}0 &
    \alpha^*\\ \alpha & 0\end{pmatrix}.
\end{equation*}
Similarly,  $i_!(x)$ corresponds to $(\cE,\widetilde{\theta})$
where:
 \begin{equation*}
    \cE = C_0(X,\xi)\underset{\rho}{\otimes}C_0(E, p^*(\Lambda E))
    \simeq C_0(E, p^*(\xi\otimes\Lambda E)) \text{ and }
    \widetilde{\theta}=
    \begin{pmatrix}
    0 & \theta^* \\ \theta & 0
    \end{pmatrix}
    \in\cL(\cE) .
 \end{equation*}
We next use the language of \cite{Blackadar, Olsen}, where
the notion of ``connection'' in the framework of Kasparov's theory
was defined. It is easy to check that
$M(\widetilde{\alpha}\hat{\otimes} 1)$ is a $0$-connection on
$\cE$ and $N(1\hat{\otimes} C)$ is a $C$-connection on $\cE$ (the
symbol $\hat{\otimes}$ denotes the graded tensor product), which
yields that:
\begin{equation*}
   \widetilde{\theta} = M(\widetilde{\alpha}\hat{\otimes} 1)+
   N(1\hat{\otimes} C)
\end{equation*}
is a $C$-connection on $\cE$. Moreover, for any $f\in C_0(X)$, we
have
\begin{equation*}
   f[\widetilde{\alpha}\hat{\otimes}1,\widetilde{\theta}]f^*
   =2M|f|^2\widetilde{\alpha}^2\hat{\otimes}1\ge 0,
\end{equation*}
which proves that $(\cE,\widetilde{\theta})$ represents the
Kasparov product of $x$ and $T$.
\end{proof}
It is known that $T$ is invertible in $KK$-theory (\cite{Ka1},
paragraph 5, theorem 8). We now give a description of its inverse
via a deformation groupoid when the bundle $E$ is the
complexification of a real euclidean bundle $N$. Hence let us
assume that $E=N\otimes\CC$ or, up to a $\CC$-linear vector bundle
isomorphism, let us assume that the bundle $E$ is the Withney sum
$N\oplus N$ of two copies of some real euclidean vector bundle $p_N :
N \to X$ with the complex structure given by $J(v,w)=(-w,v)$,
$(v,w)\in N\oplus N$. We endow the complex bundle $E$ with the induced
hermitian structure. We then define the {\sl Thom groupoid} as
follows:
\begin{equation*}
    \cI_N := E\times [0,1]
    \rightrightarrows N\times [0,1]
\end{equation*}
with structural morphism given by 
\begin{align*}
    r(v,w,0) & = s(v,w,0) = (v,0)\\
    r(v,w,t) & =(v,t),\\
    s(v,w,t) & =(w,t), \quad  t>0 \\
    (v,w,0)  \cdot (v,w',0) & = (v,w+w',0) \ \text{ and }\\
    (v,w,t)  \cdot (w,u,t)  & =(v,u,t) \quad  t>0.
\end{align*}
Thus, for $t=0$, the groupoid structure of $E$ corresponds to the
vector bundle structure given by the first projection $E=N\oplus
N\to N$ while for $t>0$ the groupoid structure of $E$ corresponds
to the pair groupoid structure in each fiber $E_x=N_x\times N_x$.

The topology of $\cI_{N}$ is inherited from the usual tangent groupoid
topology, in particular $\cI_{N}$ is a Hausdorff topological groupoid
that can be viewed as a continuous field of groupoids over $X$ with
typical fiber the tangent groupoid of the typical fiber of the vector
bundle $N \to X$.  More precisely, the topology of $\cI_{N}$ is such
that the map $E\times[0,1]\to \cI_N$ sending $(u,v,t)$ to $(u,u+tv,t)$
if $t>0$ and equal to identity if $t=0$ is a homeomorphism.

The family of Lebesgue measures on euclidean fibers $N_x$, $x\in X$,
gives rise to a continuous Haar system on $\cT_N$ that allows us to
define the $C^*$-algebra of $\cI_N$ as a continuous field of amenable
groupoids. Therefore, $\cI_N$ is amenable. We also get an element of
$KK(C^*(E),C_0(X))$, denoted by $T_{inv}$ and defined as usual by:
\begin{equation*}
    T_{inv} := [e_0]^{-1}\otimes [e_1]\otimes \cM.
\end{equation*}
Here, as before, the morphism $e_0 : C^*(\cI_N)\to C^*(\cI_N|_{t=0})=
C^*(E)$ is the evaluation at $0$, the morphism $e_1 : C^*(\cI_N)\to
C^*(\cI_N|_{t=1})$ is the evaluation at $1$, and $\cM$ is the natural
Morita equivalence between $C^*(\cI_N|_{t=1})$ and $C_0(X)$. For
instance, $\cM$ is represented by the Kasparov module $(\cH,m,0)$
where $\cH$ is the continuous field over $X$ of Hilbert spaces
$\cH_x=L^2(N_x)$, $x\in X$ and $m$ is the action of
$C^*(\cI_N|_{t=1})=C^*(N\underset{X}{\times}N)$ by compact operators
on $\cH$.

We denote $T_0=(\cE_0,\rho_0,F_0)\in KK(C_0(X),C^*(E))$ the element
corresponding to the Thom element $T$ of proposition \ref{propThom}
trough the isomorphism $C_0(E)\simeq C^*(E)$. This isomorphism is
given by the Fourier transform applied to the second factor in
$E=N\oplus N$ provided with the groupoid structure of
$\cI_N|_{t=0}$. The $C^*(E)$-Hilbert module $\cE_0=C^*(E,\Lambda E)$
is the natural completion of $C_c(E,p^*(\Lambda E))$ ($p$ is the
bundle map $E\to X$). The representation $\rho_0$ of $C_0(X)$ and the
endomorphism $F_0$ of $\cE_0$ are given by
\begin{equation*}
\begin{gathered}
    \rho_0(f)\omega(v,w)=f(x)\omega(v,w),\\
    F_0\omega(v,w) =\int_{(w',\xi)\in N_{x}
    \times  N_{x}^*}\hspace{-1cm}e^{i(w-w').\xi}
    C(v+i\xi)\omega(v,w')dw'd\xi.
\end{gathered}
\end{equation*}
In the above formulas, $f\in C_0(X)$, $\omega\in C_c(E,p^*(\Lambda
E))$ and $(v,w)\in E_x$.  We can therefore state the following result.

\begin{theorem}\label{thom-complexified-real-bundle}
The elements $T_{inv}$ and $T_0$ are inverses to each other in
$KK$--theory.
\end{theorem}

\begin{proof} We know (\cite{Ka1}, paragraph 5, theorem 8) that
$T$, hence $T_0$, is invertible so it is enough to check that
$T_0\otimes T_{inv} = 1\in  KK(C_0(X),C_0(X))$.

Since $T_{inv} := [e_0]^{-1}\otimes [e_1]\otimes \cM$ where $e_t$ are
restriction morphisms at $t=0,1$ in the groupoid $\cI_N$ we first
compute $\widetilde{T}=T_0\otimes[e_0]^{-1}$, that is, we look for
$\widetilde{T}=(\cE,\rho,F)\in KK(C_0(X),C^*(\cI_N))$ such that
\begin{equation*}
     (e_0)_*(\widetilde{T})=(\cE
   \underset{e_0}{\otimes}C^*(E),\rho,F\otimes1) = T_0
\end{equation*} 
Let $\cE=C^*(\cI_N,\Lambda E)$  be the $C^*(\cI_N)$-Hilbert
module completion of $C_c(\cI_N,(r')^*\Lambda E)$, where
$r'= p\circ \mathrm{pr}_1\circ r :\cI_N\to X$. 
Let us define a representation $\rho$ of $C_0(X)$ on $\cE$ by 
\begin{equation*}
        \rho(f)\omega(v,w,t)=f(p(v))\omega(v,w,t) \qquad 
        \text{for all } f\in C_0(X), \omega\in \cE,(v,w,t)\in\cI_N
        \end{equation*}
Let $F$ be the endomorphism of $\cE$ densely defined on  
$C_c(\cI_N,(r')^*\Lambda E)$ by 
\begin{equation*}
        F\omega(v,w,t)=\int_{(v',\xi)\in N_x\times
        N_x^*}e^{i(\frac{v-v'}{t}).\xi}
        C(v+i\xi)\omega(v',w,t)\frac{dv'}{t^{n}}d\xi,
\end{equation*}
if $t>0$ and by $F\omega(v,w,0)=F_0\omega(v,w,0)$ if $t=0$. The
integer $n$ above is the rank of the bundle $N\to X$. One can
check that the triple $(\cE,\rho,F)$ is a Kasparov
$(C_0(X),C^*(\cI_N))$--module and that under the obvious isomorphism
\begin{equation*}
  q\cE \underset{e_0}{\otimes}C^*(E)\simeq \cE_0,
\end{equation*}
$\rho$ coincides with $\rho_0$ while $F\otimes 1$ coincides with $F_0$.

\smallskip  Next, we evaluate $\widetilde{T}$ at $t=1$
and $T_1 := (e_1)_*(\widetilde{T})\in KK(C_0(X),
C^*(N\underset{X}{\times}N))$ is represented by $(\cE_1,
\rho_1,F_1)$ where $\cE_1 = C^*(N\underset{X}{\times}N,\Lambda E)$
is the $C^*(N\underset{X}{\times}N)$--Hilbert module completion of
$C_c(\cI_N|_{t=1},(p\circ r)^*\Lambda E)$ and $\rho_1,F_1$ are
given by the formulas above where $t$ is replaced by $1$.

\smallskip  Now, applying the Morita equivalence $\cM$ to $T_1$ gives:
 \begin{equation*}
   \left(\cE_1,\rho_1,F_1\right)\otimes
       \left( \cH, m, 0\right)=
    \left(\cH_{\Lambda E}, \phi, F_1\right),
 \end{equation*}
where $\cH_{\Lambda E}=(L^2(N_x,\Lambda E_x))_{x\in X}$,
$\phi$ is the obvious action of $C_0(X)$ on $\cH_{\Lambda E}$ and
$F_1$ is the same operator as above identified with a continuous
family of Fredholm operators acting on $L^2(N_x,\Lambda E_x)$:
 \begin{equation*}
   F_1\omega(x,v)=\int_{(v',\xi)\in
     N_x\times N_x^*}e^{i(v-v').\xi}
   C(v+i\xi)\omega(x,v')dv'd\xi  .
 \end{equation*}
By (\cite{CS} lemma 2.4) we know that $\left(\cH_{\Lambda E},
\phi, F_1\right)$ represents $1$ in $KK(C_0(X),C_0(X))$ (the key
point is again that the equivariant $O_n$-index of $F_1$
restricted to even forms is $1$, see also \cite{HOR69}) and the
theorem is proved.
\end{proof}

\smallskip Now let us consider the vector bundle $p: T\cV\To TM$, where $M$
is a compact manifold embedded in some $\RR^N$ and $\cV$ is the
normal bundle of the embedding. We let $q:TM \rightarrow M$ be the
canonical projection and to simplify notations, we denote again by
$p$ the bundle map $\cV\to M$ and by $\cV$ the pull-back of $\cV$
to $TM$ via $q$.

\smallskip  Using the identifications $T_xM\oplus \cV_x\simeq
T_{(x,v)}\cV$ for all $x\in M$ and $v\in \cV_x$, we get the
isomorphism of vector bundles over $TM$:
\begin{equation*}
    q^*(\cV\oplus\cV)\ni (x,X,v,w)
    \longmapsto (x,v;X+w)\in T\cV .
\end{equation*}
It follows that $T\cV$ inherits a complex structure from
$\cV\oplus\cV\simeq \cV\otimes\CC$ and we take the Atiyah-Singer
convention: via the above isomorphism, the first
parameter is real and the second is imaginary.

\smallskip The previous construction leads to the groupoid $\cI_{\cV}$ giving
the inverse of the Thom isomorphism. Actually, we slighty modify
to retain the natural groupoid structure carried by the base space
$TM$ of the vector bundle $T\cV$ (it is important in the purpose
of extending the Thom isomorphism to the singular setting). Thus,
we set:
 \begin{equation*}
  \cH_{\cV}= T\cV\times \{0\} \sqcup \pb{p}(TM)\times ]0,1]
  \rightrightarrows \cV\times [0,1] \ .
 \end{equation*}
This is the Thom groupoid defined in the section
\ref{deformation}. The groupoids $\cI_{\cV}$ and $\cH_{\cV}$ are
not isomorphic, but a Fourier transform in the fibers of $TM$
provides an isomorphism of their $C^*$-algebras:
$C^*(\cI_{\cV})\simeq C^*(\cH_{\cV})$. Moreover, this isomorphism
is compatible with the restrictions morphisms and we can rewrite
the theorem (\ref{thom-complexified-real-bundle}):

\begin{corollary}\label{inverse-thom-smooth-case}
Let $\partial_{\cH_{\cV}}=[\epsilon_0]^{-1} \otimes [\epsilon_1]$
be the $KK$-element associated with the deformation groupoid
$\cH_{\cV}$ and let $\cM$ be the natural Morita equivalence
between $C^*(\pb{p}(TM))$ and $C^*(TM)$. Then
$T_{inv}=\partial_{\cH_\cV}\otimes \cM \in KK(C^*(T\cV),C^*(TM))$
gives the inverse of the Thom isomorphism $T \in
KK(C_0(T^*M),C_0(T^*\cV)$ trough the isomorphisms $C_0(T^*M)\simeq
C^*(TM)$ and $C_0(T^*\cV)\simeq C^*(T\cV)$.
\end{corollary}

\begin{remarks}\label{bott+boundary}
1) Let us assume that $M$ is a point and $\cV=\RR^N$.
The groupoid $\cH_{\cV}$ is equal in that case to the tangent
groupoid of the manifold $\RR^N$ and the associated $KK$-element
$\partial_{\cH_\cV}\otimes \cM$ gives the Bott periodicity between the
point and $\RR^{2N}$.\\
2) Let $M_+$ be a compact manifold with boundary and $M$ the
manifold without boundary obtained by doubling $M_+$. Keeping the
notations above, let $\cV_+$ be the restriction of $\cV$ to $M_+$.
All the previous constructions applied to $M$ restrict to $M_+$ and
give the inverse $T^+_{inv}$ of the  Thom element
$T^+\in KK(C_0(T^*M_+),C_0(T^*\cV_+))$.
\end{remarks}

\smallskip  With this description of the (inverse) Thom isomorphism
in hand, the equality between the analytical and topological
indices of Atiyah and Singer \cite{AS1} follows from a commutative
diagram. Let us denote by $p_{[0,1]}$ the map $ p\times \hbox{Id}
: \cV\times[0,1]\to M\times[0,1]$. We consider the deformation
groupoid (cf. example 1. of \ref{deformation})
\begin{equation*}
  \cT^t_{\cV}= \cG^t_{\cV} \times \{0\} \sqcup
  \pb{p_{[0,1]}}(\cG^t_M) \simeq \cH_{\cV} \times \{0\} \sqcup (\cV
  \times \cV \times [0,1]) \times ]0,1] \rightrightarrows \cV \times
  [0,1] \times [0,1] \ .
 \end{equation*}
We use the obvious notation for restriction morphisms (cf. proof
of definition-proposition \ref{Tinv}) and $\cM$ for the various
(but always obvious) Morita equivalence maps. To shorten the
diagram, we set $K(G):=K_0(C^*(G))$ for all the (amenable)
groupoids met below. We have:

\begin{equation}\label{index-theorem-proof} \hspace{-0.8cm} {\small
   \begin{CD}  K(\cdot)@<=<<
  K(\cdot)@<{[ev_0]}<< K(\lbrack 0,1\rbrack) @>[ev_1]>> K(\cdot)@>=>>
  K(\cdot)\\
 @A\cM AA @A\cM AA @A\cM AA @A\cM AA @A\cM AA \\
K(\RR^N \times \RR^N) @<{j_1}<< K(\cV\times \cV) @<{[ev_0]}<<
K(\cV\times \cV \times [0,1]) @>{[ev_1]}>> K(\cV\times \cV)@>\cM>>
K(M\times M)
\\ @A{[e_1^{\RR^N}]}AA @A{[e_1^{\cV}]}AA  @AA{[q_{1,\cdot}]}A
@AA{[\pb{p}e_1^M]}A @AA{[e_1^M]}A
  \\
K(\cG^t_{\RR^N}) @<j<< K(\cG^t_{\cV})  @<{[q_{\cdot,0}]}<<
K(\cT^t_{\cV}) @>{[q_{\cdot,1}]}>>
K(\pb{p_{[0,1]}}(\cG^t_M)) @>{\cM}>> K(\cG^t_M) \\ @V{[e_0^{\RR^N}]}VV
@V{[e_0^{\cV}]}VV @VV{[q_{0,\cdot}]}V   @VV{[\pb{p}e^M_0]}V @VV{[e^M_0]}V
\\
K(T\RR^N) @<<{j_0}< K(T\cV) @<<{[\epsilon_0]}< C^*(\cH_{\cV})
@>>{[\epsilon_1]}> K(\pb{p}(T M)) @>>{\cM}> K(TM)
\end{CD} }
\end{equation}
The commutativity of this diagram is obvious. {}From the previous
remark we deduce that the  map $K(T\RR^N)\to \ZZ$ associated with
$[e^{\RR^N}_0]^{-1}\otimes [e^{\RR^N}_1]\otimes \cM$ on the left
column is equal to the Bott periodicity isomorphism $\beta$.
Thanks to the corollary (\ref{inverse-thom-smooth-case}), the map
$K(TM)\to K(T\cV)$ associated with $\cM^{-1}\otimes
\partial_{H_\cV}^{-1}=T_0$ on the bottom line is equal to the Thom
isomorphism, while $j_0$ is the usual excision map resulting from
the identification of $\cV$ with an open subset of $\RR^N$.  It
follows that composing the bottom line with the left column
produces the map: \begin{equation*}\beta\circ j_0\circ T_0:
K(TM)\to \ZZ\end{equation*} which is exaclty the Atiyah-Singer's
definition of the topological index map. We already know that the
map $K(TM)\to\ZZ$ associated with $[e^M_0]^{-1}\otimes
[e^M_1]\otimes \cM$ on the right column is the analytical index
map. Finally, the commutativity of the diagram and the fact that
the map associated with $[ev_0]^{-1}\otimes [ev_1]$ on the top
line is identity, completes our proof of the Atiyah-Singer
index theorem.

\smallskip Another proof of the usual Atiyah-Singer index theorem in the
framework of deformation groupoids can be found in
\cite{MonthubertNistor}.

\subsection{An index theorem for conical pseudomanifolds}

We define for a conical manifold a topological index and prove the
equality between the topological and analytical indices. Both
indices are straight generalisations of the ones used in the
Atiyah-Singer index theorem: indeed, if we apply our constructions
to a smooth manifold and its tangent space, we find exaclty the
classical topological and analytical indices. Thus, the egality of
indices we proove can be presented as the index theorem for
conical manifolds. Moreover, the $K$-theory of $T^{\fS}X$ is
exhausted by {\sl elliptic symbols} associated with
pseudo-differential operators in the $b$-calculus \cite{jml} and
the analytical index can be interpreted via the Poincar\'e duality
\cite{DL}, as their Fredholm index.

\smallskip  Let $X$ be a compact conical pseudomanifold embedded in
$(\RR^N)^{\fS}$ for
some $N$ and let $\cW$ be a tubular neighborhood of this embedding
as in \ref{geom}. We  first assume that $X$ has only
one singularity. We denote by
\begin{equation*}
    \cH_{\cW}=T^{\fS} \cW \times\{0\} \sqcup
    \pb{\pi}(T^{\fS}X)\times]0,1]
    \rightrightarrows \cW^{\circ}\times [0,1]
\end{equation*}
the Thom groupoid associated with $\pi : \cW \to X$ and by
\begin{equation*}
    \cH_+ = T \cW_+\times \{0\} \sqcup \pb{\pi_+}(TX_+)\times]0,1]
    \rightrightarrows \cW_+ \times [0,1]
\end{equation*}
the Thom groupoid associated with $\pi_+:\cW_+\to
X_+$. Here $\cW_+ = \cW \setminus O_{\cW}=\{(z,V)\in \cW \
\vert \ \rho(z)\geq 1\}$ and $X_+=X\setminus O_X=\{z\in X \ \vert
\ \rho(z)\geq 1\}$, where $\rho$ is in both case the defining
function of the singularity. We denote by $T_{inv}$ and $T_{inv}^+$ the respective
inverse-Thom elements. Recall (cf. prop. \ref{shortsequence}) that
we have the two following short exact sequences comming from
inclusion and restriction morphisms:
\begin{equation*}
\begin{CD}
    0 \longrightarrow K(\cK(L^2(O_{ \cW })))
    @>{i_{\cV}}_*>> K(C^*(T^\fS  \cW ))
    @>{r_{\cW}}_*>> K(C^*(T \cW_+)) \longrightarrow 0 \\
    0 \longrightarrow K(\cK(L^2(O_X))) @>>> K(C^*(T^\fS X))
    @>>> K(C^*(T X_+)) \longrightarrow 0
   \end{CD} \
\end{equation*}
\begin{dfprop}\label{invertibility-Tc} The following diagram commutes:
\begin{equation}\label{thom-thom+}
    \begin{CD} 0 \longrightarrow K(\cK(L^2(O_{ \cW })))
    @>{i_{ \cW }}_*>> K(C^*(T^\fS  \cW ))
    @>{r_{ \cW }}_*>> K(C^*(T \cW_+)) \longrightarrow 0 \\
    @V\cM VV @V \cdot \otimes T_{inv} VV
    @V \cdot \otimes T_{inv}^+ VV \\
    0 \longrightarrow K(\cK(L^2(O_X))) @>>> K(C^*(T^\fS X))
    @>>> K(C^*(T X_+)) \longrightarrow 0
    \end{CD}
\end{equation}
where $\cM$ is the natural Morita equivalence map. In
particular, the map:
\begin{equation}\label{thomc}
\cdot \otimes T_{inv} : K(C^*(T^\fS  \cW ))\To
  K(C^*(T^\fS X))
\end{equation}
is an isomorphism. Its inverse is denoted by $T$ and called the Thom
isomorphism.
\end{dfprop}
\begin{proof}
Let us note again by $\pi$ the (smooth) vector bundle map
$\cW^{\circ} \to \smx$ and consider the following diagram:
\begin{equation}\label{thom-thom+-proof}
   \begin{CD}
      0 \to C^*(O_{ \cW }\times O_{ \cW }) @>>> C^*(T^\fS  \cW )
      @>r_+>> C^*(T \cW_+) \to 0 \\
      @Aev_0AA    @Aev_0AA    @Aev_0AA  \\
      0 \to C^*(O_{ \cW } \times O_{ \cW } \times
      \lbrack 0,1 \rbrack)
      @>>> C^*(\cH_\cW) @>r_+>> C^*(\cH_+) \to 0 \\
      @Vev_1VV  @Vev_1VV   @Vev_1VV  \\
      0 \to C^*(\pb{\pi}(O_X\times O_X))
      @>>> C^*(\pb{\pi}(T^\fS X))
      @>r_+>> C^*(\pb{\pi}(TX_+)) \to 0
   \end{CD}
\end{equation}
where the (Lie) groupoid isomorphism $\pb{\pi}(O_X\times
O_X)\simeq O_{ \cW }\times O_{ \cW }$ has been used. Applying the
$K$ functor and Morita equivalence maps to the bottom line to get
rid of the pull back $\pb{\pi}$ and using the fact that the long
exact sequences in $K$-theory associated to the top and bottom
lines split in short exact sequences, give the diagram
(\ref{thom-thom+}). Since $\cM$ and $T_{inv}^+$ are isomorphisms,
the same is true for $T_{inv}$.
\end{proof}
\begin{remarks}
1)  When $X$ has several singular points, the invertibility of
$.\otimes T_{inv}$ remains true. This can be checked thanks to a
recursive process on the number $k$ of singular points. First choose a
singular point $s\in \fS$ and call again $s$ its image in $\cW$ by
the embedding $\pi$. Denote by 
\begin{equation*}
    \cH_{s,+} = T \cW_{s,+}\times \{0\} \sqcup \pb{\pi_{s,+}}(T^{\fS} X_{s,+})\times]0,1]
    \rightrightarrows \cW_{s,+} \times [0,1]
\end{equation*}
the Thom groupoid associated with $\pi_{s,+}:\cW_{s,+}\to
X_{s,+}$. Recall that $\cW_{s,+} = \{(z,V)\in \cW \
\vert \ \rho_s(z)\geq 1\}$ and $X_{s,+}=X\setminus O_s=\{z\in X \ \vert
\ \rho_s(z)\geq 1\}$, where $\rho_s$ is in both case the defining
function associated to $s$. We denote by $T_{inv}^{s,+}$ the corresponding
inverse-Thom element. The same proof as before gives that the map:
\begin{equation*}
\cdot \otimes T_{inv} : K(C^*(T^\fS  \cW ))\To
  K(C^*(T^\fS X))
\end{equation*}
is an isomorphism as soon as \begin{equation*}
\cdot \otimes T_{inv}^{s,+} : K(C^*(T^\fS  \cW _{s,+}))\To
  K(C^*(T^\fS X_{s,+}))
\end{equation*} is. But now $X_{s,+}$ has $k-1$ singular points. 

\noindent 2) The Thom map we define extends the usual one: this is
exactly what is said by the commutativity of the diagram
(\ref{thom-thom+}).
\end{remarks}

\smallskip Let us recall that we started with an embedding of $X$ into
$(\RR^N)^{\fS}$ which is $\RR^N$
with $k$ singular points where $k$ is the cardinal of $\fS$. The $\fS$-tangent space $T^\fS \cW $ of $
\cW $ is obviously isomorphic to an open subgroupoid of the
$\fS$-tangent space $T^\fS(\RR^N)^{\fS}$. Thus we get an excision
homomorphism:
\begin{equation*}
    j : C^*(T^\fS  \cW ) \To C^*(T^\fS(\RR^N)^{\fS}).
\end{equation*}
There is a natural identification of the $K$-theory group $K(T^\fS
{(\RR^N)^{\fS}})$ with $\ZZ$, analog to the one given by Bott
periodicity in the case of $T\RR^N=\RR^{2N}$ coming from its
tangent groupoid (cf. remark \ref{bott+boundary}):
  \begin{equation*}
    \partial_{(\RR^N)^{\fS}}= [e_0]^{-1}\otimes [e_1]\otimes \cM :
    K(T^\fS (\RR^N)^{\fS}) \To \ZZ
  \end{equation*}
  \begin{equation*}
    K(C^*(T^\fS (\RR^N)^{\fS}))\overset{[e_0]}{\longleftarrow}
    K(C^*(\cG^t_{(\RR^N)^{\fS}})\overset{[e_1]}{\longrightarrow}
    K(\cK(L^2(\RR^N)))\overset{\cM}{\rightarrow} K(\cdot)\simeq\ZZ
  \end{equation*}
We are now in position to extend the Atiyah-Singer topological
index to conical pseudomanifolds:

\begin{definition} The topological index of the conical
  pseudomanifold $X$  is defined by
\begin{equation*}
   \ind^X_t =
   \partial_{(\RR^N)^{\fS}}\circ [j] \circ T
\end{equation*}
\end{definition}
Moreover, we obtain the following extension of the
Atiyah-Singer Index theorem

\begin{theorem} If $X$ is a pseudomanifold with conical
singularities then
\begin{equation*}
    \ind^X_a=\ind^X_t
\end{equation*}
\end{theorem}

\begin{proof}
The proof is similar to our proof of the Atiyah-Singer index
theorem. Indeed, let us write down the analog of the diagram
(\ref{index-theorem-proof}) for the singular manifold $X$:
 
\begin{equation}\label{conical-index-theorem-proof} \hspace{-0.8cm} {\small
   \begin{CD}  K(\cdot)@<=<<
  K(\cdot)@<{[ev_0]}<< K(\lbrack 0,1\rbrack) @>[ev_1]>> K(\cdot)@>=>>
  K(\cdot)\\
 @A\cM AA @A\cM AA @A\cM AA @A\cM AA @A\cM AA \\
K(\RR^{2N}) @<{j_1}<< K(\cW^2) @<{[ev_0]}<<
K(\cW^2\times [0,1]) @>{[ev_1]}>> K(\cW^2)@>\cM>>
K((X^o)^2)
\\ @A{[e_1^{(\RR^N)^c}]}AA @A{[e_1^{\cW}]}AA  @AA{[q_{1,\cdot}]}A
@AA{[\pb{\pi}e_1^X]}A @AA{[e_1^X]}A
  \\
K(\cG^t_{(\RR^N)^{\fS}}) @<j<< K(\cG^t_{\cW})  @<{[q_{\cdot,0}]}<<
K(\cT^t_{\cW}) @>{[q_{\cdot,1}]}>>
K(\pb{\pi_{[0,1]}}(\cG^t_X)) @>{\cM}>> K(\cG^t_X) \\ @V{[e_0^{(\RR^N)^c}]}VV
@V{[e_0^{\cW}]}VV @VV{[q_{0,\cdot}]}V   @VV{[\pb{\pi}e^X_0]}V @VV{[e^X_0]}V
\\
K(T^\fS (\RR^N)^{\fS}) @<<{j_0}< K(T^\fS\cW) @<<{[\epsilon_0]}< C^*(\cH_{\cW})
@>>{[\epsilon_1]}> K(\pb{\pi}(T^\fS X)) @>>{\cM}> K(T^\fS X)
\end{CD} }
\end{equation}
This diagram involves various deformation groupoids associated to $X$
and its embedding into $(\RR^N)^{\fS}$. The commutativity is obvious since
everything comes from morphisms of algebras or from explicit Morita
equivalences. As before, the convention $K(G)=K_0(C^*(G))$ is used to
shorten the diagram and intuitive notations are chosen to name the
various restriction morphisms. Starting from the bottom right corner
and following the right column gives the analytical index
map. Starting from the bottom right corner and following the bottom
line and next the left column gives the topological index map.
\end{proof}

%%%%%%%%%%%%%%%

\subsection{Signification of the index map}
In the sequel we suppose that $X$ has only one singularity.\\
In \cite{DL}, a Poincar\'e duality in bivariant $K$-theory
between $C(X)$ and $C^*(T^\fS X)$ is proved. Taking the Kasparov product with 
the dual-Dirac element involved in
this duality provides an isomorphism:  
\begin{equation}\label{abstractsymbolmap}
    K_0(X)\overset{\Sigma_X}{\longrightarrow} K^0(T^\fS X)
\end{equation}
When a $K$-homology class  of $X$ and a $K$-theory class of $T^\fS X$
coincide trough this isomorphism, we say that they are Poincar\'e
dual. 

If $p : X\to \cdot$ is the trivial map, then (\ref{abstractsymbolmap}) satisfies \cite{jml}:
\begin{equation}\label{Ind-PD-maptoapoint}
  \ind^X_a \circ \Sigma_X = p_* : K_0(X)\longrightarrow  \ZZ \simeq K_0(\cdot)
\end{equation}
Remember that cycles of $K_0(Y)$, for a compact Hausdorff space $Y$, are
given by triples $(H,\pi,F)$ where $H=H_+\oplus H_-$ is a $\ZZ_2$-graded Hilbert
space, $\pi$ a degree $0$ homomorphism of $C(Y)$ into the algebra of
bounded operators on $H$ and $F=\begin{pmatrix}0 & F_-\\ F_+ & 0\end{pmatrix}$ 
a bounded operator on $H$ of degree $1$ such that $F^2-1$ and $[\pi,F]$ are compact. 
Since:
\begin{equation*}
  p_*(H,\pi,F) = \text{Fredholm-Index}(P),
\end{equation*}
the equality (\ref{Ind-PD-maptoapoint}) implies that $\ind^X_a$ produces indices of
Fredholm operators. To make things more concrete and see what Fredholm
operators come into the play, one needs to compute explicitely
(\ref{abstractsymbolmap}), or mimeting the case of smooth manifolds, interpret it as a {\sl
  symbol map} associating {\sl $K$-theory classes of the tangent
  space} to {\sl elliptic pseudodifferential
  operators}. 
This is done in full details, and summarized below, for the $0$-order
case in \cite{jml}. We give also an account of the unbounded case,
necessary to compare our index with the ones computed in
\cite{cheeg,bs,chou,lesch}.  

\subsubsection{$K$-homology of the conical pseudomanifold and elliptic
  operators}
Let $\Psi_b^*$ be the algebra of the $b$-calculus \cite{meaps} on
$\overline{X^\circ}$ (the obvious compactification of $X^\circ$ into a
manifold with boundary). A $b$-pseudodifferential operator $P$ is said to
be {\sl fully elliptic} if its principal symbol $\sigma_{int}(P)$ , regarded as an
ordinary pseudodifferential operator on $X^\circ$ is invertible and
the indicial family $(\hat{P}(\tau))_{\tau\in\RR}$  is everywhere
invertible \cite{meaps} (that is, for all $\tau\in\RR$, the pseudodifferential
operator $\hat{P}(\tau)$ on $L=\partial\overline{X^\circ}$ is
invertible). A {\sl full parametrix} of $P$ is then another
$b$-operator $Q$ such that $PQ$ and $QP$ are equal to $1$ modulo a
negative order $b$-operator with vanishing indicial family. When $P$
is a zero order fully elliptic $b$-operator, it is Fredholm on the
Hilbert space $L^{2,b}:=L^2(X^\circ, d\mu_b)$ for the natural measure
$d\mu_b=\frac{dh}{h}dy$ coming with an exact $b$-metric \cite{meaps}, 
and we get a canonically defined $K$-homology class of $X$:
\begin{equation}\label{K-hom-via-boperators}
  [P] := [(H^b,\pi,\mathbf{P})]\in K_0(X)
\end{equation}
where $\mathbf{P}=\begin{pmatrix} 0 & Q \\ P & 0\end{pmatrix}$, 
$H^b=L^{2,b}\oplus L^{2,b}$ and $\pi : C(X)\to \cB(H^b)$ is the homomorphism given by pointwise
multiplication (for all $f\in C^\infty_c(X^\circ)\oplus\CC$, $[\pi(f),\mathbf{P}]$ has negative order and
vanishing indicial family, thus it is a compact operator on $H^b$
\cite{meaps} ; since 
$C^\infty_c(X^\circ)\oplus\CC$ is dense in $C(X)$, it proves that the
commutators $[\pi,\mathbf{P}]$ are compact and $[P]$ is well defined). 

\subsubsection{$K$-theory of the noncommutative tangent space and symbols}
To compute the Poincar\'e dual of $K$-classes given by (\ref{K-hom-via-boperators}),
one uses a slightly different, but $KK$-equivalent, definition of $T^\fS X$: 
\begin{equation}
  \label{TqX}
   T^\fS X := T]0,1[\times L\times L\sqcup TX_+\rightrightarrows X^{o}
\end{equation}
The $KK$-equivalence between both definitions is explicit (\cite{jml})
and allows us to translate all the previous
constructions to this
variant of the tangent space. 

Roughly speaking, a noncommutative
symbol on the pseudomanifold $X$ is a pseudodifferential operator, in
the groupoid sense (\cite{Monthubert,NWX}), on $T^\fS X$. For
technical reasons, one asks to these objects to be smooth up to $h=0$,
in other words we define the algebra of noncommutative symbols as:
\begin{equation}
  \label{noncomm-symb-def}
  S^*(X)= \Psi^*(\overline{T^\fS X})
\end{equation}
where $\overline{T^\fS X}=\{0\}\times\RR\times L\times L \cup T^\fS X$
and the letter $\Psi$ is reserved for the space of pseudodifferential
operators on the indicated groupoid. See \cite{jml} for the precise
assumptions on the Schwartz kernels of the operators in (\ref{noncomm-symb-def}).
Considering the closed saturated subspace $L=\partial
\overline{X^\circ}$ of the space of units $\overline{X^\circ}$ of
$\overline{T^\fS X}$, we get a restriction homomorphism: 
\begin{equation}
  \label{restriction-hom-psidiff}
    S^*(X)= \Psi^*(\overline{T^\fS X}) \overset{\rho}{\longrightarrow}
    \Psi^*(\RR\times L\times L)\simeq \Psi^*_{susp}(L)
\end{equation}
where $\Psi^*_{susp}(L)$ denotes the space of suspended
pseudodifferential operators of R. Melrose \cite{Mel2}. 
A noncommutative symbol $a\in S^m(X)$ is {\sl fully elliptic} if there
exists $b\in S^{-m}(X)$ such that $ab$ and $ba$ are equal to $1$
modulo $S^{-1}(X)\cap \ker\rho=:\cJ$. Fully elliptic symbols $a\in S^0(X)$
give canonically $K$-classes of the tangent space $T^\fS X$: 
\begin{equation}
  \label{K-class-symbol-order0}
   [a] = [\cE,\mathbf{a}]\in KK(\CC,C^*(T^\fS X))=K^0(T^\fS X)
\end{equation}
 where $\mathbf{a}=\begin{pmatrix} 0 & b\\a & 0\end{pmatrix}$, $b$ any
 inverse of $a$ modulo $\cJ$ and $\cE=C^*(T^\fS X)\oplus C^*(T^\fS
 X)$.

\subsubsection{$\ind^\fS_a$ as a Fredholm index}
The main result of \cite{jml} is:
\begin{theorem}\label{main-jml}
  There exists a surjective linear map $\sigma_X : \Psi_b^*\to S^*$ such that:
  \begin{itemize}
  \item $P\in\Psi^*_b$ is fully elliptic if and only if $\sigma_X(P)$
    is fully elliptic,
  \item For all zero order fully elliptic operator $P$,
    \begin{equation}
      \label{concrete-PD}
         \Sigma_X([P])=[\sigma_X(P)]
    \end{equation}
  \end{itemize}
\end{theorem}
See also \cite{LMN1} for a thorough study of the property of full
ellipticity in $b$-calculus in the framework of groupoids. 
\begin{remarks}
\item Allowing vectors bundles $E$ over $\overline{X^\circ}$ and
defining the algebra of $b$-operators $\Psi^*_b(E)$ and the
algebra of noncommutative symbols
$S^*(X,E)$ accordingly, we get a full description of
$K_0(X)$ in terms of $b$-operators and of $K^0(T^\fS X)$ in terms of
noncommutative symbols. This is also proved in \cite{jml}. Thus, for any $x\in K^0(T^\fS X)$, we have:
 $$  \ind^X_a(x) = \hbox{Fredholm-index}(P_x) $$
 where $P_x$ is any $b$-operator such that $[\sigma_X(P_x)]=x$.
 \item The reader should not be
surprised by our definition of (noncommutative) symbols: if $V$ is
smooth manifold, the algebra of ordinary symbols is isomorphic to the
algebra of pseudodifferential operators, for a suitable choice of
regularizing operators imposed by the use of the Fourier transform, on
the {\sl groupoid} $TV$. 
 
  \end{remarks}

\subsubsection{The unbounded case and geometric operators}
The symbol map $\sigma_X$ constructed in \cite{jml} makes sense on
differential $b$-operators. It turns out that natural geometric
operators on $X$, when provided with a conical metric, can be written
as $b$-differential operators with singular coefficients at $h=0$, or,
in the terminology of \cite{lesch}, Fuchs type operators. We
explain in this paragraph how to relate the analysis of these
operators (\cite{cheeg,bs,chou,lesch}) to our $K$-theoritic
constructions, for the case of a Dirac operator on $X$ even
dimensional with one conical point $s$. 

Let $g$ be a riemmannian metric on $X^\circ$ which is conical on
$O_X=]0,1[\times L$, that is: $g=dh^2+h^2g_L$. We assume to simplify
the computations that the riemannian metric $g_L$ on $L$ is
independant of $h$ when $h\le 1$. We denote by $d\hbox{vol}_{X}$ the
corresponding volume form.  

Let $\cE=\cE_+\oplus\cE_-$ be a Clifford module and $c$ the
corresponding Clifford multiplication (\cite{bgv}). 

If $X^\circ$ has a spin structure, then there exists a
($\ZZ_2$-graded) vector bundle $\cW$ such that $\cE\simeq \cW\otimes \cS$ 
where $\cS=\cS_+\oplus\cS_-$ is the spinor bundle.  

In the case of a spin structure, using the canonical metric structure
and Clifford connection $\nabla^{\cS}$ of 
the spinor bundle, and using on $\cW$ a metric structure of product
type on $O_x$ and a compatible connection $\nabla^\cW$, we get on
$\cE$ a metric structure, such that $\cE_+\perp\cE_-$ and
$c(v)^*=-c(v)$ for unitary tangent vectors $v\in TX^\circ$, and a Clifford
connection $\nabla^{\cE}$ such that the corresponding Dirac operator
$D$ is symmetric, when considered as an unbounded operator on
$L^2(\cE)$ with domain the space $C^\infty_c(\cE)$ of compactly
supported sections. Recall that $D$ is defined locally by the formula:
 $$
   s\in C^\infty(\cE),\qquad Ds = \sum_{i=1}^n c(e_i)\nabla^{\cE}_{e_i}s,
 $$
where $(e_1,\ldots,e_n)$ is a local basis of $TX^\circ$. 

If $X^\circ$ has no spin
structure, then  the isomorphism $\cE\simeq \cW\otimes \cS$ remains
true locally. Thus, one can still construct locally on $\cE$ metrics
and connections with the previous properties and then patch them with
a partition of unity $(U_i,\phi_i)$ on $X$ (that is $(U_i)$ is a
finite open covering of $X^\circ$ by open charts, $\phi_i\in
C^\infty(U_i)\cap C_c(U_i\cup \{s\})$ and $\sum_i\phi(x)=1,\ \forall x\in X^\circ$). The 
resulting Dirac operator is again symmetric, and all subsequent
computations are exactly the same with or without spin structure. 

Although the Hilbert space $L^2(\cE)$, whose scalar product is given by
\begin{equation}
  \label{eq:conical-hilbert}
 (s,t)=\int_{X^\circ} \left(s(x),t(x)\right)_{\cE_x}\ d\hbox{vol}_X  
\end{equation}
is the most natural Hilbert
space with respect to the given geometric data, computations are easier
with $H_b(E)$ which is defined as, if
$\pi : X^\circ\to X_+$ denotes the obvious retraction map and $\widetilde{\cE}=\cE|_L$, the
completion of $C^\infty_c(\pi^*(\cE|_{X_+}))$ for the scalar product:
\begin{equation}
  \label{eq:b-hilbert}
  (s,t)_b=\int_{X_+} \left(s(x),t(x)\right)_{\cE_x}\ d\hbox{vol}_X + 
    \int_{x=(h,y)\in O_X} \!\!\!\!\!\!
    \left(s(x),t(x)\right)_{\widetilde{\cE}}\ \frac{dh}{h}d\hbox{vol}_Y(y)
\end{equation}
One can choose an isometry $U:H_b(E)\to L^2(E)$ such that 
$U :C^\infty(\pi^*(\cE|_{X_+}))\to C^\infty(\cE)$ is equal to identity on
the complement of some open neighborhood of $\overline{O_X}$ and given
on $O_X$ by
 $$ U(s)(h,y) = h^{-\frac{n}{2}}\theta^{\cE}_{(1,y)\to (h,y)}s(h,y) $$
where $\theta^{\cE}_{(1,y)\to (h,y)}:\cE_{(1,y)}\to\cE_{(h,y)}$ is the
parallel transport associated with the connection and the canonical
identification 
$C^\infty(\pi^*(\cE|_{X_+})|_{O_X})\simeq C^\infty(]0,1[,C^\infty(\widetilde{\cE}))$
have been used.  

Then, a straight computation shows that (\cite{bc,chou,jml0}) the
following holds for sections $s$ supported on $O_X$:
\begin{equation}
  \label{eq:Dirac-Fuchs}
  U^{-1}DU s = c(e_1)\cdot\frac{\partial
    s}{\partial h}+
  \frac{1}{h}\left(\widetilde{D} -\frac{e_1}{2}\right)s
\end{equation}
where $(e_1=\frac{\partial}{\partial h},e_2,\ldots,e_n)$ is a local
orthonormal basis in $TO_X$ and $\widetilde{D}$ is the differential
operator on $L$, acting on the sections of $\widetilde{\cE}$, given by 
$\widetilde{D} u = \sum_{i=2}^n c(\widetilde{e_i})\nabla^{\widetilde{\cE}}_{\widetilde{e_i}}u$,
where $\widetilde{e_i}(y) =e_i(1,y)$ and $\nabla^{\widetilde{\cE}}$ is the
connection on $\widetilde{\cE}$ induced by $\nabla^{\cE}$. Moreover we
have $\cE_-=c(e_1)\cdot\cE_+$, and  the operator $U^{-1}DU$  
is given by the matrix , in the decomposition 
$\cE=\cE_+\oplus c(e_1)\cdot\cE_+\simeq \cE_+^2$,  
\begin{equation}
  \label{eq:Fuchs}
 \begin{pmatrix}
    0 & -\frac{\partial}{\partial h} + \frac{1}{h}(S + \frac{1}{2}) \\
    \frac{\partial}{\partial h} + \frac{1}{h}(S - \frac{1}{2}) & 0
 \end{pmatrix}
\end{equation}
where $S=\sum_{i=2}^n c(\widetilde{e_i}\widetilde{e_1})\nabla^{\widetilde{\cE}}_{\widetilde{e_i}}$
is again a symmetric Dirac operator on the Clifford module
$\widetilde{\cE_+}$ over $L$. 

It is of course equivalent to study $D$ on $L^2(E)$ or $T=U^{-1}DU$ on
$H_b(E)$. If we used $dh$ instead of $\frac{dh}{h}$ in
(\ref{eq:b-hilbert}), which leads to the Hilbert space used in
\cite{bs,lesch} to study Fuchs type operators, $S$ would appear without the extra terms
$\pm \frac{1}{2}$ in (\ref{eq:Fuchs}).

The deformation process of \cite{jml} used to associate noncommutative
symbols to $b$-pseudodifferential operators, can be applied to $T$ and
gives a family $(T_t)_{0\le t\le 1}$ where 
$t>0,\ T_t\in \frac{1}{h}.\Psi^1_b(\cE)$ and 
 $T_0\in \frac{1}{h}.S^{1}(X,\cE)$ have the following expression on
 $O_X$: 
\begin{equation}
  \label{eq:Tt}
 t>0,\quad  T_t = \begin{pmatrix}
    0 & -t\frac{\partial}{\partial h} + \frac{1}{h}(S + \frac{t}{2}) \\
    t\frac{\partial}{\partial h} + \frac{1}{h}(S - \frac{t}{2}) & 0
 \end{pmatrix}   
\end{equation}
and 
\begin{equation}
  \label{eq:T0}
  \sigma_X(T):=T_0= \frac{1}{h}\begin{pmatrix}
    0 & -\frac{\partial}{\partial \lambda} + S \\
    \frac{\partial}{\partial \lambda} + S  & 0
 \end{pmatrix}
\end{equation}
Observe that the family $D_t:=UT_tU^{-1}$ coincides for $t>0$ with the
one given by the deformation of the conical metric
$\frac{dh^2}{t^2}+h^2g_L$ (\cite{bc}).

The natural questions are then: does the noncommutative symbol $T_0$ give
canonically a $K$-theory element of $T^\fS X$ ? does the operator $T$ give
canonically a $K$-homology class on $X$ ? Are the corresponding
classes Poincar\'e dual ? 

The answer to the first two questions is negative in general, but becomes
affirmative under some conditions on the spectrum of $S$ and in
that case the answer to the last question is affirmative too. Let us
explain these phenomena. 

Firstly, the noncommutative symbol $a:=hT_0$ is fully elliptic if and only
if 
\begin{equation}
  \label{eq:0notinthespectrum}
  0 \not\in\hbox{spec} S
\end{equation}
We assume in the sequel that
this condition is satisfied. Using the ellipticity of $a$ as a
pseudodifferential operator on $\overline{T^\fS X}$ and the
invertibility of $a|_{h=0}$, we can proove thanks to \cite{Vas} that 
$(1+\sigma_X(T)^2)^{-1}\in h^2\cdot S^{-2}(X,\cE)\subset \cK(C^*(T^\fS X,\cE))$. 
This implies that the closure of $\sigma_X(T)$ as an unbounded
operator on the $C^*(T^\fS X)$-Hilbert module $C^*(T^\fS X,\cE)$ with
domain $C^\infty(T^\fS X,\cE)$ is selfadoint, regular and provides  an unbounded 
$(\CC,C^*(T^\fS X))$-Kasparov bimodule (\cite{bj,Vas}). We thus get here a
well defined, canonical, element $[\sigma_X(T)]\in K^0(T^\fS X)$. 

We turn back now to the operators $T_t$, $1\ge t>0$. It is well known that
they always have a selfadjoint extension, not unique in general
(\cite{bs,chou,lesch}). Adpating for instance the computations of
\cite{bs} to our particular choice of Hilbert space $H_b(\cE)$, we see
that $T_t$, $t>0$, with domain $C^\infty_c(E)$ is essentially self-adjoint if
and only if $\hbox{spec}(S)\cap ]-\frac{t}{2},\frac{t}{2}[=\emptyset$.
Otherwise, any choice of an orthogonal decomposition of:
\begin{equation}
  \label{eq:space-of-choices}
   W_t=\bigoplus_{-t/2< u< t/2}\CC.e_u
\end{equation}
where the $e_u$'s describe an orthonormal system of eigenvectors of
$S$, allows to define a self-adjoint extension of $T$ (\cite{bs}). 

Thus, for $\alpha$ small enough, thanks to the assumption
(\ref{eq:0notinthespectrum}), $T_\alpha$ is essentially
self-adjoint. It is also Fredholm by (\cite{bs}), so $T_{\alpha}$ gives
an unbounded $(C(X),\CC)$-Kasparov bimodule, in other words a
$K$-homology class $[T_\alpha]\in K_0(X)$, and we have: 
 $$ 
    \Sigma_X([T_\alpha])=[\sigma_X(T)]
 $$
To check this, one shows that the Woronowicz transform
$q(T_\alpha)=T_\alpha.(1+T_\alpha^2)^{-1/2}$ (\cite{Vas,bj})
of $T_\alpha$ can be represented in $K$-homology by a zero order
$b$-operator with noncommutative symbol equal to the Woronowicz
transform $q(\sigma_X(T))=\sigma_X(T).(1+\sigma_X(T)^2)^{-1/2}$ of $\sigma_X(T)$, and then the theorem
\ref{main-jml} applies. 

In particular: 
 $$ 
  \dim\ker(T_\alpha)_+-\dim\ker(T_\alpha)_-=\ind^X_a([\sigma_X(T)])
 $$
In general, $T=T_1$ has several selfadjoint extensions, but using the
splitness of 
$$ 
 0 \longrightarrow \CC \longrightarrow K_0(X)\longrightarrow
 K_0(X^\circ)\longrightarrow 0
$$
one shows that two given selfadjoint extensions of $T$ give the same
$K$-homology  class if and only if their Fredholm index is the same. 
Thus a selfadjoint extension $T_Z$, given by a choice of a
decomposition $Z\oplus Z^\perp$ of (\ref{eq:space-of-choices}),
produces the  same $K$-homology class as $T_\alpha$ (and then, is Poincar\'e dual to its
noncommutative symbol) if and only if $2\dim Z=\dim W_1$. 

Let us say a word about the case $0\in\hbox{spec}S$. 
For small $t$, the selfadjoint extensions of $T_t$ are classified by the
orthogonal decompositions of $\ker S$. There is a priori no canonical
choice. On the other hand,  the noncommutative symbol $\sigma_X(T)$ is not
fully elliptic. We conjecture that the selfadjoint
extensions of $\sigma_X(T)$, as an unbounded operator on the
Hilbert module $C^*(T^\fS X,\cE)$, are again classified 
by the orthogonal decomposition of $\ker S$ and give unbounded
Kasparov modules which are in one-to-one correspondance, via Poincar\'e
duality, with the selfadjoint extensions of $T_t$. 
 
\bibliographystyle{plain}

%\bibliographystyle{plain}
%\bibliography{Itop}
\end{document}